\documentclass[11pt]{article}

\usepackage{amsmath}
\usepackage{amssymb}
\usepackage{latexsym}

\newtheorem{assumption}{Assumption}

\def\qed{ \ \vrule width.2cm height.2cm depth0cm\smallskip}

\newenvironment{proof}{\noindent {\bf Proof.\/}}{$\qed$\vskip 0.1in}

\newcommand{\la}{\langle}
\newcommand{\ra}{\rangle}

\newcommand{\eps}{\varepsilon}

\newcommand{\ba}{\begin{array}}
\newcommand{\ea}{\end{array}}
\newcommand{\be}{\begin{equation}}
\newcommand{\ee}{\end{equation}}
\newcommand{\bea}{\begin{eqnarray}}
\newcommand{\eea}{\end{eqnarray}}
\newcommand{\beaa}{\begin{eqnarray*}}
\newcommand{\eeaa}{\end{eqnarray*}}

\def\dbE{\mathbb{E}}
\def\dbF{\mathbb{F}}

\def\dbL{\mathbb{L}}

\def\dbP{\mathbb{P}}

\def\dbR{\mathbb{R}}
\def\dbS{\mathbb{S}}

\def\dbU{\mathbb{U}}
\def\dbV{\mathbb{V}}

\def \dbf{{\mathbf{d}}}

\def\a{\alpha}
\def\b{\beta}
\def\g{\gamma}
\def\d{\delta}
\def\e{\varepsilon}

\def\si{\sigma}
\def\t{\tau}
\def\f{\varphi}
\def\th{\theta}
\def\o{\omega}

\def\G{\Gamma}
\def\D{\Delta}

\def\L{\Lambda}

\def\O{\Omega}

\def\cA{{\cal A}}
\def\cB{{\cal B}}
\def\cC{{\cal C}}

\def\cF{{\cal F}}

\def\cH{{\cal H}}

\def\cO{{\cal O}}
\def\cP{{\cal P}}

\def\cT{{\cal T}}
\def\cU{{\cal U}}
\def\cV{{\cal V}}

\def\cY{{\cal Y}}
\def\cZ{{\cal Z}}
\def\no{\noindent}

\def\ms{\medskip}

\def\q{\quad}
\def\qq{\qquad}

\def\pa{\partial}
\def\cd{\cdot}
\def\cds{\cdots}

\def\ch{\textsc{h}}
\def\tr{\hbox{\rm tr}}

\def\qed{ \hfill \vrule width.25cm height.25cm depth0cm\smallskip}

\newcommand{\basa}{\begin{assumption}}
\newcommand{\easa}{\end{assumption}}

\newcommand{\bas}{\begin{assum}}
\newcommand{\eas}{\end{assum}}

\def\pa{\partial}

 \def\cd{\cdot}
\def\cds{\cdots}

\def\as{\hbox{\rm -a.s.{ }}}

\def\tr{\hbox{\rm tr$\,$}}

\def\wh{\widehat}
\def\dis{\displaystyle}

\def\1{\mathbf{1}}

\def\:{\!:\!}
\def\reff#1{{\rm(\ref{#1})}}
\def \proof{{\noindent \it Proof.\quad}}

\newcommand{\ol}{\overline}
\newcommand{\ul}{\underline}

at 9pt

\begin{document}

\newtheorem{thm}{Theorem}[section]
\newtheorem{lem}[thm]{Lemma}
\newtheorem{cor}[thm]{Corollary}
\newtheorem{prop}[thm]{Proposition}
\newtheorem{rem}[thm]{Remark}
\newtheorem{eg}[thm]{Example}
\newtheorem{defn}[thm]{Definition}
\newtheorem{assum}[thm]{Assumption}

\renewcommand {\theequation}{\arabic{section}.\arabic{equation}}
\def\thesection{\arabic{section}}

\title{\bf Two Person Zero-sum Game in Weak Formulation and Path Dependent Bellman-Isaacs Equation}

\author{Triet {\sc Pham}\footnote{University of Southern California, Department of Mathematics, trietpha@usc.edu. Research supported by USC Graduate School Dissertation Completion Fellowship.}
        \and Jianfeng {\sc Zhang}\footnote{University of Southern California, Department of Mathematics, jianfenz@usc.edu. Research supported in part by NSF grant DMS 10-08873.  Part of the research was done while this author was visiting Shandong University, whose hospitality is greatly appreciated.} \footnote{The authors would like to thank Rainer Buckdahn and Lihe Wang for very helpful discussion. In particular we thank Buckdahn for providing us the counterexample Example \ref{eg-Buckdahn}, and Wang for showing us the proof of Lemma \ref{lem-Wang}.} 
}\maketitle

\begin{abstract}
In this paper we study  a two person zero sum stochastic differential game in weak formulation. Unlike standard literature which uses strategy type of controls, the weak formulation allows us to consider  the game with control against control. We shall prove the existence of game value under natural conditions. Another main feature of the paper is that we allow for non-Markovian structure, and thus the game value is a random process. We characterize the value process as the unique viscosity solution of the corresponding path dependent Bellman-Isaacs equation, a notion recently introduced by Ekren, Keller, Touzi and Zhang \cite{EKTZ} and Ekren, Touzi and Zhang \cite{ETZ1, ETZ2, ETZ3}. 
\end{abstract}



\section{Introduction}
\label{sect-Introduction}
\setcounter{equation}{0}

Since the seminal paper  Fleming and Souganidis \cite{FS}, two person zero sum stochastic differential games have been studied extensively in the literature, see e.g. \cite{BY}, \cite{BLions}, \cite{BCQ},  \cite{BHL}, \cite{BL},   \cite{CR}, \cite{EH}, \cite{ES}, \cite{HL}, \cite{HW}, \cite{KSud}, \cite{MY}, \cite{Swiech}, to mention a few. There are typically two approaches. One is to use the viscosity theory, namely to show that the value function of the game is the unique viscosity solution of the associated Bellman-Isaacs equation, and the other is to use the Backward SDE approach, which characterizes the value process as the solution to a related  BSDE.

To be precise, let $u$ and $v$ denote the controls of the two players,  $B$ a Brownian motion,  $X^{S,u,v}$ the controlled state process in the strong formulation:
\bea
\label{Xstrong}
X^{S,u,v}_t = x + \int_0^t b(s,  u_s, v_s)  ds + \int_0^t \si(s,  u_s, v_s)  dB_s,
\eea
and $J(u,v)$ the corresponding value (utility or cost) which is determined by $X^{S,u,v}$, $B$, and $(u,v)$. The lower and upper values of the game are defined as:
\beaa
\underline V_0 := \sup_{u\in \cU}\inf_{v\in \cV} J(u,v),\q  \overline V_0 := \inf_{v\in \cV}\sup_{u\in \cU} J(u,v),
\eeaa
where $\cU$ and $\cV$ are appropriate sets of admissible controls. It is clear that $\underline V_0\le \overline V_0$. Two central problems in the game literature are:

{\it (i) When does the game value exists, namely $V_0:=\underline V_0 =  \overline V_0$?

(ii) Given the existence of the game value, is there a saddle point? That is, we want to find $(u^*, v^*)\in \cU\times \cV$ such that $V_0 = J(u^*,v^*) = \inf_{v\in\cV} J(u^*,v) = \sup_{u\in\cU} J(u, v^*)$.}

\ms

\no However, even under reasonable assumptions, the game value may not exist. We shall provide a counterexample, see Example \ref{eg-Buckdahn} below, which is due to Buckdahn. 

To overcome the difficulty,    Fleming and Souganidis \cite{FS} introduced strategy types of controls:
\beaa
\underline V'_0 := \sup_{\a\in \cA}\inf_{v\in \cV} J(\a(v),v),\q  \overline V_0' := \inf_{\b\in \cB}\sup_{u\in \cU} J(u,\b(u)),
\eeaa
Here $\a : \cV\to \cU$ and $\b: \cU \to \cV$ are so called strategies  and $\cA$, $\cB$ are appropriate sets of admissible strategies. Under the Isaacs condition and assuming the comparison principle for the viscosity solution of the corresponding Bellman-Isaacs equation holds, \cite{FS} showed that $\underline V_0' = \overline V_0'$. This work has been extended by many authors in various aspects. In particular, Buckdahn and Li \cite{BL} defined $J(u,v)$ via Backward SDEs, and very recently Bayraktar and Yao \cite{BY} used doubly reflected BSDEs. The main drawback of this approach, however, is that the two players have non-symmetric information, and for $\underline V'_0$ and $ \overline V_0'$, the roles of two players are switched. Consequently, it is less convenient to study the saddle point in this setting. 
 
We propose to attack the problem in weak formulation, which  is more convenient for proving the Dynamic Programming Principle. Note that in \reff{Xstrong} the controls $(u, v)$ actually mean $u(B_\cd), v(B_\cd)$. Our weak formulation is equivalent to the following feedback type of controls:
 \bea
\label{Xweak}
X^{W, u,v}_t = x + \int_0^t b(s,  u_s(X^{W,u,v}_\cd), v_s(X^{W,u,v}_\cd))  ds  + \int_0^t \si(s,   u_s(X^{W,u,v}_\cd), v_s( X^{W,u,v}_\cd))  dB_s,
\eea
Here $X^{W,u,v}_\cd$ denotes the path of $X^{W,u,v}$ and the superscript $^W$ stands for weak formulation. Under natural assumptions, we show that the game value does exist. The advantage of the weak formulation setting is that we are using control against control, thus one can define the saddle point naturally. When there is only drift control, namely $\si$ is independent of $(u,v)$, one can prove the existence of saddle point under mild conditions. However, when there is diffusion control, the problem is much more involved. We shall obtain some approximate saddle point.

We remark that, when there is only drift control, the weak formulation has already been used in the literature, see Bensoussan and Lions \cite{BLions}  for Markovian case and Hamadene and Lepetier \cite{HL} for non-Markovian case.  The former one relies on PDE arguments and the latter one uses Backward SDEs. The advantage in this case is that one can easily obtain the weak solution of SDE \reff{Xweak} by applying the Girsanov Theorem. Our general case with diffusion control has different nature. Roughly speaking, the drift control is associated with semi-linear PDEs, while the diffusion control is associate with fully nonlinear PDEs. We also note that, in a Markovian model but also with optimal stopping problem, Karatzas and Sudderth \cite{KSud} studied the game problem with diffusion control in weak formulation, under certain strong conditions.

Another main feature of our paper is that we study the game in non-Markovian framework, or say in a path dependent manner. The standard approach in the literature, e.g. \cite{FS} and \cite{BL}, is to prove that the lower value and the upper value are a viscosity solution (or viscosity semi-solution) of the corresponding Bellman-Isaacs equation, then by  assuming the comparison principle for the viscosity solution of the PDE, one obtains the existence of the game value. These works rely on the PDE arguments and thus works only in Markovian setting. In a series of papers, Ekren, Keller, Touzi and Zhang \cite{EKTZ} and Ekren, Touzi and Zhang \cite{ETZ1, ETZ2, ETZ3} introduced a notion of viscosity solution for the so called path dependent PDEs and established its wellposedness. This enables us to extend the above approach to  path dependent setting. Indeed, based on the dynamic programming principle we establish, we show that the lower value and the upper value of the game are viscosity solutions of the corresponding path dependent  Bellman-Isaacs equations. Then, under the Isaacs condition and assuming the uniqueness of viscosity solutions, we characterize the game value as the unique viscosity solution of  the path dependent  Bellman-Isaacs equation.

Finally we remark that, due to weak formulation with diffusion control, this paper is by nature closely related to the  second order BSDEs (2BSDEs, for short) introduced by Cheridito, Soner, Touzi and Vicoir \cite{CSTV} and Soner, Touzi and Zhang \cite{STZ1, STZ2}, and the $G$-expectation introduced by Peng \cite{Peng-G}. While more involved here, our arguments for Dynamic Programming Principle follow the idea in \cite{STZ1, STZ2} and Peng \cite{Peng-g}. However, $G$-expectations and 2BSDEs  involve only stochastic optimization and thus  the generator is convex in terms of the hessian. Consequently, the dynamic value process is a supermartingale under each associated probability measure.  For our game problem, the Bellman-Isaacs equation is non-convex, and the value process is not a supermartingale anymore. Under additional technical conditions, we conjecture that our value process will be a semi-martingale. This requires to develop a semi-martingale theory under nonlinear expectation and to generalize the 2BSDE theory  to non-convex generators. We established some norm estimates for semi-martingales in another paper Pham and Zhang \cite{PZ} and will leave the general 2BSDE theory for future research.

The rest of the paper is organized as follows. In Section \ref{Sect - Prelim} we present some preliminaries.  The game problem is introduced in Section \ref{sect-game}.  In Sections \ref{sect-DPP} and \ref{sect-Viscsol}  we  prove the dynamic programming principle and the viscosity property, respectively. In Section \ref{sect-comparison} we study the comparison principle for PPDEs and in Section \ref{sect-saddle} we investigate approximate saddle points. Finally some technical proofs are presented Appendix.  

\section{Preliminaries}
\label{Sect - Prelim}
\subsection{The canonical space}

Let $\O:= \big\{\o\in C([0,T], \dbR^d): \o_0={\bf 0}\big\}$, the set of continuous paths starting from the origin, $B$ the canonical process, $\dbF$ the filtration generated by $B$, 
$\dbP_0$ the Wiener measure,  and $\L := [0,T]\times \O$. Here and in the sequel, for notational simplicity we use ${\bf 0}$ to denote vectors or matrices with appropriate dimensions whose components are all equal to $0$. Let $\dbS^d$ denote the set of $d\times d$ matrices, $\dbS^d_{\ge {\bf 0}} := \{\si \in \dbS^d: \si\ge {\bf 0}\}$, 
and
 \beaa
 &x \cd x' := \sum_{i=1}^d x_i x'_i 
 ~~\mbox{for any}~~x, x' \in \dbR^d,
 ~~\g : \g' := \mbox{Trace}[\g\g']
 ~~\mbox{for any}~~\g, \g'\in \dbS^d.
 &
 \eeaa
 We define a norm on $\O$ and a metric on $\L$ as follows: for any $(t, \o), ( t', \o') \in\L$,
\bea\label{norm}
 \|\o\|_{t} 
 := 
 \sup_{0\le s\le t} |\o_s|,
 \q  
 \dbf_\infty\big((t, \o),( t', \o')\big) 
 := 
 |t-t'| + \big\|\o_{.\wedge t} - \o'_{.\wedge t'}\big\|_T.
 \eea
Then $(\O, \|\cd\|_{T})$ and $(\L, \dbf_\infty)$ are complete metric spaces.  
\begin{defn}
\label{defn-spaceC0} Let  $Y: \L\to \dbR$ be  an $\dbF$-progressively measurable process.
\\
{\rm (i)}  We say $Y\in \dbL^\infty(\L)$ if $Y$ is bounded. 
\\
{\rm (ii)} We say $Y\in C^0(\L)$ (resp. $UC(\L)$) if $Y$ is continuous (resp. uniformly continuous) in $(t,\o)$. Moreover, we denote $C^0_b(\L) := C^0(\L) \cap \dbL^\infty(\L)$ and $UC_b(\L) := UC(\L) \cap \dbL^\infty(\L)$.
\\
{\rm (iii)} We say $Y \in {\underline \cU}$  if $Y$ is bounded from above, upper semi-continuous (u.s.c. for short)  from right in $t$, and there exists a modulus of continuity function $\rho$ such that for $(t,\o), (t',\o')\in\L$:
 \bea\label{USC}
 Y(t,\o) - Y(t', \o') 
 \le 
 \rho\big(\dbf_\infty((t,\o),(t',\o'))\big)
 ~\mbox{whenever}~t\le t',
 \eea 
and we say $Y \in {\overline \cU}$ if $-Y \in {\underline \cU}$.
\end{defn}

It is clear that ${\underline \cU} \cap {\overline \cU} = UC_b(\L)$. Moreover, we denote by $\dbL^\infty(\L, \dbR^d)$  the space of $\dbR^d$-valued processes whose components are in $\dbL^\infty(\L)$, and define other similar notations  in the same spirit.

\ms

We next introduce the shifted spaces. Let  $0\le s\le t\le T$.

-  Let $\O^t:= \big\{\o\in C([t,T], \dbR^d): \o_t ={\bf 0}\big\}$ be the shifted canonical space; $B^{t}$ the shifted canonical process on
$\O^t$;   $\dbF^{t}$ the shifted filtration generated by $B^{t}$, $\dbP^t_0$ the Wiener measure on $\O^t$, and $\L^t := [t,T]\times \O^t$.

- Define $\|\cd\|_s$ on $\O^t$, $\dbf_\infty$ on $\L^t\times\L^t$,  and $C^0(\L^t)$ etc. 
in the spirit of (\ref{norm}) and Definition \ref{defn-spaceC0}.

- For  $\o\in \O^s$ and $\o'\in \O^t$, define the concatenation path $\o\otimes_{t} \o'\in \O^s$ by:
\beaa
(\o\otimes_t \o') (r) := \o_r\1_{[s,t)}(r) + (\o_{t} + \hat\o'_r)\1_{[t, T]}(r),
&\mbox{for all}&
r\in [s,T].
\eeaa

- Let $s\in[0,T)$ and $\o\in \O^s$. For an $\cF^{s}_{T}$-measurable random variable $\xi$, an $\dbF^{s}$-progressively measurable
process $X$ on $\O^s$, and $t\in(s,T]$, define the shifted $\cF^{t}_{T}$-measurable random variable $\xi^{t,\o}$  and  $\dbF^{t}$-progressively measurable
process $X^{t,\o}$ on $\O^t$ by:
\beaa
\xi^{t, \o}(\o') :=\xi(\o\otimes_t \o'), \q X^{t, \o}(\o') := X(\o\otimes_t \o'),
&\mbox{for all}&
\o'\in\O^t.
\eeaa

It is clear that, for any $(t,\o) \in \L$ and any $Y\in \dbL^\infty(\L)$, we have $Y^{t,\o} \in \dbL^\infty(\L^t)$. Similarly the property holds for other spaces defined in Definition \ref{defn-spaceC0}.

\subsection{Probability measures}
In this subsection we introduce the probability measures on $\O^t$ in different formulations. First, let $\si \in \dbL^\infty(\L, \dbS^d_{\ge {\bf 0}})$, $b\in \dbL^\infty(\L, \dbR^d)$. Define
\bea
\label{strong}
\dbP^{S,\si, b} := \dbP_0\circ (X^{S, \si, b})^{-1} &\mbox{where}& X^{S, \si, b}_t := \int_0^t b_s ds +  \int_0^t \si_s dB_s,\q \dbP_0\mbox{-a.s.}
\eea
Here the superscript $^S$ stands for strong formulation.    We next introduce the corresponding weak formulation. We denote  a probability measure $\dbP$ on $\O$ as $\dbP^{W, \si, b}$ if 
\bea
\label{weak}
M^{b}_t := B_t - \int_0^t b_s ds ~\mbox{is a $\dbP$-martingale and}~\la M^{b}\ra_t = \int_0^t \si_s^2 ds~~\dbP\mbox{-a.s.}
\eea
Here the quadratic variation $\la M^{b}\ra$ is under $\dbP$. We remark that $\dbP^{W, \si, b} := \dbP_0\circ (X^{W, \si, b})^{-1}$, where $X^{W, \si, b}$  is a weak solution of the following SDE  (with random measurable coefficients): 
\bea
\label{weaksolution}
X^{W, \si, b}_t :=   \int_0^t b_s(X^{W, \si, b}_\cd) ds + \int_0^t \si_s(X^{W, \si, b}_\cd) dB_s,\q \dbP_0\mbox{-a.s.}
\eea 
In other words, we are considering feedback type of controls.

In this paper we shall use the weak formulation, which is more convenient for proving Dynamical Programming Principle. We note that, for arbitrarily given $(\si, b)$, the SDE \reff{weaksolution} may not have a weak solution, namely there is no $\dbP$ such that $\dbP = \dbP^{W,\si, b}$. Let 
\bea
\label{cAWS}
\left.\ba{lll}
\overline \Xi^W := \Big\{(\si, b)\in \dbL^\infty(\L, \dbS^d_{\ge {\bf 0}}) \times \dbL^\infty(\L, \dbR^d):  \mbox{ SDE \reff{weaksolution}  has a unique weak solution}\Big\};\\
\overline \Xi^S := \Big\{(\si, b)\in \dbL^\infty(\L, \dbS^d_{\ge {\bf 0}}) \times \dbL^\infty(\L, \dbR^d):  \mbox{ SDE \reff{weaksolution}  has a unique strong solution}\Big\}.
\ea\right.
\eea
%
%
For probability measures on the shifted space $\O^t$, we define $\dbP^{S, t, \si, b}$, $\dbP^{W, t, \si, b}$, and $\overline\Xi^{W,t}$, $\overline \Xi^{S,t}$, etc.  similarly.

We next introduce the regular conditional probability distribution (r.c.p.d for short) due to  Stroock and Varadhan \cite{SV}. We shall follow the presentation in Soner, Touzi and Zhang \cite{STZ1}. Let $\dbP$ be an arbitrary probability measure on $\O$ and $\t$ be an $\dbF$- stopping time. The r.c.p.d.  $\{\dbP^{\t,\o}, \o\in\O\}$ satisfies: 
\begin{itemize}
\item For each $\o$, $\dbP^{\t,\o}$ is a probability measure on $\cF^{\t(\o)}_T$;
\item  For every bounded $\cF_T$-measurable random variable $\xi$:
\bea
\label{rcpd}
\dbE^{\dbP} \big[\xi | \cF_\t\big] (\o) = \dbE^{\dbP^{\t,\o}}\big[\xi^{\t(\o), \o}\big],\q  \dbP \as
\eea
\end{itemize}

The following simple lemma will be important for the proof of Dynamic Programming Principle in Section \ref{sect-DPP} below. Its proof is postponed to Appendix. 

\begin{lem}
\label{lem-paste}
Let $(\si, b)\in\overline \Xi^{S}$ (resp. $\overline \Xi^W$), $t\in [0,T]$, $\{E_i, 1\le i\le n\}\subset \cF_t$ be a partition of $\O$, and $(\si^i, b^i)\in \overline \Xi^{S,t}  (\mbox{resp.} ~\overline \Xi^{W,t})$. Define
\beaa
\bar\si (\o) := \si(\o) \1_{[0,t)} + \sum_{i=1}^n \si^i(\o^t) \1_{E^i} \1_{[t,T]},\q \bar b (\o) := b(\o) \1_{[0,t)} + \sum_{i=1}^n b^i(\o^t) \1_{E^i} \1_{[t,T]}.
\eeaa
Then $(\bar\si, \bar b) \in\overline\Xi^S$ (resp. $\overline \Xi^W$), and, for $i=1,\cds, n$,
\beaa
\dbP^{\bar\si, \bar b} = \dbP^{\si, b} ~\mbox{on}~ \cF_t &\mbox{and}& \big(\dbP^{\bar\si, \bar b} \big)^{t, \o} = \dbP^{t, \si^i, b^i} ~\mbox{for}~\dbP^{\si, b}\mbox{-a.e.}~\o\in E_i.
\eeaa
\end{lem}

\subsection{Viscosity solutions of path dependent PDEs}

Our notion of viscosity solutions of Path Dependent PDEs (PPDEs for short) is introduced by Ekren, Keller, Touzi and Zhang \cite{EKTZ} for semilinear PPDE and Ekren, Touzi and Zhang \cite{ETZ1, ETZ2} for fully nonlinear PPDE. We follow the presentation in \cite{ETZ1, ETZ2} here.

For any constant $L>0$, denote
\bea
\label{cPL}
\cP_L^t := \big\{ \dbP^{W,t,\si, b}:  |b|\le L, {\bf 0}\le \si \le \sqrt{2L} I_d\big\}, &\mbox{and}& \cP^t_\infty := \cup_{L>0} \cP_L^t.
\eea
We remark that in $\cP^t_L$ we do not require the uniqueness of weak solution. 

Let $Y\in C^0(\L)$. For $t\in [0,T)$,  we define the right time-derivative, if it exists, as in Dupire \cite{Dupire} and Cont and Fournie \cite{CF}:
 \bea
 \label{pat}
 \pa_t Y(t,\o)
 &:=&
 \lim_{h\downarrow 0}
 \frac{1}{h}\big[ Y\big(t+h, \o_{\cd\wedge t}\big)
                  - Y\big(t, \o\big)\big].
 \eea
For the final time $T$, we define, whenever the following limit exists:
\bea
\pa_t Y(T,\o)
&:=&
 \lim_{t\uparrow T}\pa_t  Y(t,\o).
\eea

\begin{defn}
\label{defn-spaceC12}  {\rm (i)} We say $Y\in C^{1,2}(\L)$ if $Y\in C^0(\L)$, $\pa_t Y \in C^0(\L)$, and there exist $\pa_{\o} Y \in C^0(\L, \dbR^d)$, $\pa^2_{\o\o} Y\in C^0(\L, \dbS^d)$ such that, for any $(s,\o)\in [0,T)\times \O$ and any $\dbP\in \cP^s_\infty$, $Y^{s,\o}$ is a local  $\dbP$-semimartingale and it holds:
\bea
\label{Ito}
d Y^{s,\o}_t = (\pa_t Y_t)^{s,\o}  dt+ (\pa_{\o} Y_t)^{s,\o} \cd d B^s_t + \frac12(\pa^2_{\o\o} Y_t)^{s,\o} : d \la B^s\ra_t,~~\dbP\mbox{-a.s.}
\eea
{\rm (ii)} We say $Y\in C^{1,2}_b(\L)$ if $Y\in \mbox{UC}_b(\L)$, $\pa_t Y\in C^0_b(\L)$, and the above  $\pa_{\o} Y$ and $\pa^2_{\o\o} Y$ exist and are in $C^0_b(\L, \dbR^d)$ and $C^0_b(\L, \dbS^d)$, respectively.
\end{defn}

Next, let  $\cT$  denote the set of  $\dbF$-stopping times,  and $\cH\subset\cT$ the subset of those hitting times $\ch$  taking the following form: for  some open and convex set $O  \subset \dbR^d$ containing ${\bf 0}$ and  some $0< t_0\le T$,
\bea
\label{cT}
\ch := \inf\{t: B_t \in O^c\} \wedge t_0 = \inf\{ t: d(\o_t, O^c) = 0\} \wedge t_0.
\eea

We may define $C^{1,2}(\L^t)$, $C^{1,2}_b(\L^t)$, $\cT^t$, and $\cH^t$ similarly. It is clear that, for any $(t,\o)$ and $Y\in C^{1,2}(\L)$ (resp. $Y\in C^{1,2}_b(\L)$), we have $Y^{t,\o}\in C^{1,2}(\L^t)$(resp. $Y^{t,\o}\in C^{1,2}_b(\L)$), and for any $\ch \in \cH$ such that $\ch(\o) > t$, we have $\ch^{t,\o} \in \cH^t$.

For any $L>0$, $(t,\o) \in \L$ with $t<T$, and  $\dbF$-adapted process $Y$, define
 \bea
 \label{def-cA}
 \left.\ba{lll}
 \underline\cA^{\!L} Y(t,\o)
 &:=&
\dis  \Big\{\f\in C^{1,2}_b(\L^{\!t}): ~\mbox{for some}~\ch\in \cH^t,\\
  &&\dis     (\f-Y^{t,\o})(t,{\bf 0})
       =
      \inf_{\t\in \cT^t} \inf_{\dbP \in \cP_L^t} \dbE^\dbP\big[(\f-Y^{t,\o})_{\t\wedge\ch}
                       \big]
 \Big\},
 \\
 \overline\cA^{\!L} Y(t,\o)
 &:=&
\dis  \Big\{\f\in C^{1,2}_b(\L^{\!t}): ~\mbox{for some}~\ch\in \cH^t,\\
  &&\dis     (\f-Y^{t,\o})(t,{\bf 0})
       =
      \sup_{\t\in \cT^t} \sup_{\dbP \in \cP_L^t} \dbE^\dbP\big[(\f-Y^{t,\o})_{\t\wedge\ch}
                       \big]
 \Big\}.
 \ea\right.
 \eea 
We are now ready to introduce the viscosity solution of PPDEs. Consider the following PPDE with generator $G$:
\bea
\label{PPDE}
-\pa_t Y_t - G(t,\o, Y_t, \pa_\o Y_t, \pa^2_{\o\o} Y_t) =0.
\eea

\begin{defn}
\label{defn-viscosity}
\no {\rm (i)} Let $L>0$. We say $Y\in {\underline \cU}$ (resp. ${\overline \cU}$) is a viscosity $L$-subsolution (resp. $L$-supersolution) of PPDE (\ref{PPDE})  if,  for any $(t,\o)\in [0, T)\times \O$ and any $\f \in \underline\cA^{L} Y(t,\o)$ (resp. $\f \in \overline\cA^{L}Y (t,\o)$):
 \beaa
 \Big(-\pa_t \f  - G^{t,\o}(.,  Y^{t,\o},\pa_{\o} \f,\pa^2_{\o\o}\f) \Big)(t,{\bf 0})
 &\le  ~~(\mbox{resp.} \ge)&  0.
 \eeaa

\no {\rm (ii)} We say $Y\in {\underline \cU}$ (resp. ${\overline \cU}$)  is a viscosity subsolution (resp. supersolution) of PPDE (\ref{PPDE}) if  $Y$ is viscosity $L$-subsolution (resp. $L$-supersolution) of PPDE (\ref{PPDE}) for some $L>0$.

\no {\rm (iii)} We say $Y\in UC_b(\L)$ is a viscosity solution of PPDE (\ref{PPDE})   if it is both a viscosity subsolution and a viscosity  supersolution.
\end{defn}

\begin{rem}
\label{rem-viscosity1}
{\rm For $0 < L_1 < L_2$, obviously $\cP_{L_1}^t \subseteq \cP_{L_2}^t$ and $\underline\cA^{\!L_2} Y(t,\o) \subseteq \underline\cA^{\!L_1} Y(t,\o)$. Then one can easily check that a viscosity $L_1$-subsolution must be a viscosity $L_2$-subsolution. Consequently, $Y$ is a viscosity solution of PPDE (\ref{PPDE}) iff there exists an $L \geq 1$ such that for all $\tilde{L} \geq L$, $Y$ is a viscosity $\tilde{L}$-subsolution. However, we require the same $L$ for all $(t,\o)$. A similar statement holds for the viscosity supersolution.}
\end{rem}

\begin{rem}
\label{rem-viscosity2}
{\rm (i)  In the Markovian case, namely $Y(t,\o) = Y(t,\o_t)$ and $G= g(t,\o_t,y,z,\g)$, our definition of viscosity solution is stronger than the standard viscosity solution in PDE literature. That is, $Y$ is a viscosity solution in our sense implies it is a viscosity solution in the standard sense as in Crandall, Ishii, and Lions \cite{CIL}. 

(ii) The state space $\L$ of PPDEs is not locally compact, and thus the standard arguments by using Ishii's lemma do not work in path dependent case. The main idea of \cite{EKTZ, ETZ1, ETZ2} is to transform the definition to an optimal stopping problem in \reff{cA}, which helps to obtain the comparison and hence the uniqueness of  viscosity solutions. 
\qed}
\end{rem}

\section{The zero-sum game}
\label{sect-game}
\setcounter{equation}{0}

\subsection{The admissible controls}
Let $\dbU$ and $\dbV$ be two Borel measurable spaces equipped with some topology. From now on we shall fix two $\dbF$-progressively measurable mapping:
\beaa
\si: [0, T]\times \dbU \times \dbV \to \dbS^d_{\ge {\bf 0}},&& b: [0, T] \times \dbU \times \dbV \to \dbR^d.
\eeaa

We shall always assume
\begin{assum}
\label{assum-sib}
$\si$ and $b$ are bounded by a constant $C_0$. 
%
\end{assum}

For $t\in [0, T]$, let $\overline \cU_t$ (resp. $\ol \cV_t$) denote the set of $\dbU$-valued (resp. $\dbV$-valued), $\dbF^t$-progressively measurable processes $u$ (resp. $v$) on $\L^t$. Throughout the paper, when $(u, v) \in \overline \cU_t \times \overline \cV_t$ is given, for any process $\f$ on $\L^t$ with appropriate dimension, we denote
\bea
\label{whf}
\wh \f_s := \wh \f^{t,u,v}_s := \f_s \si(s, u_s, v_s).
\eea 
Define
\bea
\label{cA}
\left.\ba{c}\Xi_t := \Big\{(u,v)\in \overline \cU_t\times \overline \cV_t:  (\si(\cd, u, v), \widehat b(\cd, u,v)) \in \overline\cA^{S,t}\Big\} \\
X^{t, u,v} := X^{W, t, \si(u, v), \widehat b(u,v)},\q  \dbP^{t, u, v} := \dbP^{W,t, \si(u, v), \widehat b(u,v)},\q \mbox{for}~(u,v) \in \Xi_t.
\ea\right.
\eea
We note that, from now on, the $\si, b$ in previous section will actually be $\si(t, u_t, v_t)$ and $\widehat b(t, u_t, v_t)$ for some $(u, v) \in \cA_0$.  In particular,  for the convenience of studying the BSDE later, we are considering SDE in the form
\bea
\label{Xweak2}
X^{t, u,v}_s = \int_t^s \si(r, u_r(X^{t,u,v}_\cd), v_r(X^{t,u,v}_\cd))\Big[ dB^t_r +   b(r, u_r(X^{t,u,v}_\cd), v_r(X^{t,u,v}_\cd))dr\Big],\q \dbP^t_0\mbox{-a.s.}
\eea
Moreover,   one can easily check that there exists a $\dbP^{t, u, v}$-Brownian motion $W^{t, u, v}$ such that
\bea
\label{Wuv}
d B^t_s = \si(s, u_s, v_s) \Big[dW^{t, u, v}_s + b(s, u_s, v_s) ds\Big],~~\dbP^{t, u,v}\mbox{-a.s.}
\eea

To formulate the game problem, we shall restrict the controls to subsets  $\cU_t \subset \overline \cU_t$ and $\cV_t \subset \overline \cV_t$ whose elements $u$ and $v$ take the following form:
\bea
\label{cUV}
&& \qq\qq u = \sum_{i=0}^{m-1} \sum_{j=1}^{n_i} u_{ij} \1_{E^i_j} \1_{[t_i, t_{i+1})},\q v = \sum_{i=0}^{m-1} \sum_{j=1}^{n_i} v_{ij} \1_{E^i_j} \1_{[t_i, t_{i+1})},~\mbox{where} \\
&& t=t_0 <\cds<t_m=T, ~\{E^i_j\}_{1\le j\le n_i}\subset \cF^t_{t_i} ~\mbox{is a partition, and $u_{ij}$, $v_{ij}$ are constants.} \nonumber
\eea
It is clear that, for $u\in \overline \cU_t$,
\bea
\label{cU}
u\in \cU_t  &\mbox{if and only if}& u ~\mbox{takes finitely many values.}
\eea
We have the following simple lemma whose proof is provided in Appendix for completeness. 
\begin{lem}
\label{lem-cUV} 
(i) $\cU_0$ is closed under pasting. That is, for $u\in \cU_0$, $t\in [0, T]$, $u^i \in \cU_t$, $i=1,\cds, n$, and disjoint $\{E_i, i=1,\cds, n\}\subset \cF_t$, the following $\ol u$ is also in $\cU_0$:
\beaa
\ol u := u \1_{[0,t)} + \Big[\sum_{i=1}^n u^i(\o^t) \1_{E_i} + u \1_{\cap_{i=1}^n E_i^c} \Big]\1_{[t, T]}.
\eeaa

(ii) Under Assumption \ref{assum-sib}, it holds $\cU_t\times  \cV_t \subset \Xi_t$.
\end{lem}
\proof In light of \reff{cU}, (i) is obvious. To see (ii), we notice that any pair of constant processes $(u, v)$ is obviously in $\Xi_t$. Then (ii) follows from repeated use of Lemma \ref{lem-paste}.
\qed

\subsection{The Backward SDEs}
Let  $f(t,\o,y,z,u,v): \L \times \dbR \times \dbR^d  \times \dbU \times \dbV \rightarrow \dbR$ be an $\dbF$-progressively measurable nonlinear generator. Throughout the paper, we shall assume
\begin{assum}
\label{assum-F}
(i) $f(t,\o,0,{\bf 0},u,v)$ is bounded by a constant $C_0$, and uniformly continuous in $(t,\o)$ with a modulus of continuity function $\rho_0$.

(ii) $f$ is uniformly Lipschitz in $(y,z)$ with a Lipschitz constant $L_0 $.
\end{assum}

Now for any $(t,\o)\in \L$,  $(u,v) \in \cU_{t} \times \cV_{t}$,  $\t\in\cT^t$, and $\cF^{t}_{\t}$-measurable  terminal condition $\eta$,  recall the notation \reff{whf} and consider  the following BSDE on $[t, \t]$: 
\bea
\label{BSDEuv}
\cY_s = \eta + \int_s^{\t} f^{t,\o}(r,B^{t}_\cd, \cY_r, \wh\cZ_r, u_r, v_r) dr - \int_s^{\t} \cZ_r dB^t_r, \; \dbP^{t,u,v} \as
\eea
We have the following simple lemma whose proof is presented in Appendix for completeness.

\begin{lem}
\label{lem-BSDE}
Let Assumptions \ref{assum-sib} and \ref{assum-F} (ii) hold, and  
\beaa
I_0^2(t,\o, u, v) := \dbE^{\dbP^{t, u, v}}\Big[|\eta|^2 + \int_{t}^{\t} |f^{t,\o}(s,B^t_\cd,0,{\bf 0},u_s,v_s)|^2ds\Big] <\infty.
\eeaa
 Then BSDE \reff{BSDEuv} has a unique solution, denoted as  $\big(\cY^{t,\o,u,v}[\t, \eta], \cZ^{t,\o,u,v}[\t, \eta]\big)$, and there exists a constant $C$, depending only on  $C_0$, $L_0$, $T$,  and the dimension $d$,  such that
\bea
\label{BSDEest}
\dbE^{\dbP^{t,u,v}}\Big[\sup_{t\le s\le \t} |\cY^{t,\o, u, v}_s[\t, \eta]|^2 + \int_{t}^{\t} |\wh \cZ^{t,\o, u, v}_s[\t, \eta]|^2ds\Big]&\le& CI_0^2(t,\o, u, v).
\eea

Moreover,  if $\t \le t+\d$, then
\bea
\label{BSDEestp}
\dbE^{\dbP^{t,u,v}}\Big[\sup_{t\le s\le \t} |\cY^{t,\o, u, v}_s[\t, \eta]|\Big] &\le& C\Big(\dbE^{\dbP^{t, u, v}}[|\eta|^2]\Big)^{1\over 2} + C \d^{1\over 2}I_0(t,\o, u, v).
\eea
\end{lem}
Throughout the paper, we shall use the generic constant $C$ which depends only on $C_0$, $L_0$, $T$,  and the dimension $d$, and may vary from line to line.  
\subsection{The value processes}
 We now fix an $\cF_T$-measurable terminal condition $\xi$ and assume throughout the paper:
\begin{assum}
\label{assum-xi}
$\xi$ is bounded by a constant $C_0$, and is uniformly continuous in $\o$ with a modulus of continuity function $\rho_0$.
\end{assum}

We now define the lower value and upper value of the game as follows:
\bea
\label{semi-value}
\underline Y(t,\o) := \sup_{u \in \cU_t} \inf_{v\in \cV_t} \cY^{t,\o,u,v}_t[T,\xi^{t,\o}];&& \overline Y(t,\o) :=  \inf_{v\in \cV_t} \sup_{u \in \cU_t} \cY^{t,\o,u,v}_t[T,\xi^{t,\o}].
\eea
As a direct consequence of Lemma \ref{lem-BSDE}, we have
\bea
\label{bound}
-C \;\le\; \underline Y &\le& \overline Y\;\le\; C.
\eea
When there is no confusion, we will simplify the notations:
\bea
\label{BSDEnotation}
 \big(\cY^{t,\o,u,v}, \cZ^{t,\o,u,v}\big):= \big(\cY^{t,\o,u,v}[T, \xi^{t,\o}], \cZ^{t,\o,u,v}[T, \xi^{t,\o}]\big).
\eea
Our goal of this paper is to show, under certain additional assumptions, that $\underline Y = \overline Y$ and it is the unique viscosity solution of certain PPDE. See Theorem \ref{thm-viscosity} below.

\begin{rem}
\label{rem-weak}
{\rm (i) In this paper we restrict our controls to $\cU_t\times \cV_t \subset \Xi_t$. We note that in general $\overline\cU_t\times \overline\cV_t$ is not in $\Xi_t$. We may study the following problem though:
\beaa
\underline Y'(t,\o) := \sup_{u \in \overline\cU_t} \inf_{v\in \overline\cV_t(u)} \cY^{t,\o,u,v}_t, && \overline Y'(t,\o) :=  \inf_{v\in \overline\cV_t} \sup_{u \in \overline\cU_t(v)} \cY^{t,\o,u,v}_t,
\eeaa
where
\beaa
\overline\cV_t(u) :=\Big\{v\in \overline\cV_t: (u, v) \in \Xi_t\Big\},\q \overline\cU_t(v) :=\Big\{u\in \overline\cU_t: (u, v) \in \Xi_t\Big\},
\eeaa
and we take the convention that, for the empty set $\phi$,  $\sup_{\phi} [\cd] = \infty$ and $\inf_{\phi}[\cd] = -\infty$. However, we will be able to prove only partial Dynamic Programming Principle in this formulation.



(ii) Another important constraint we impose is that $\si$ and $b$ are independent of $\o$.  When $\si$ and $b$ are random,  given $(t,\o)\in \L$, the solution $X^{t,\o, u, v}$ of SDE \reff{Xweak2} and its distribution $\dbP^{t,\o, u, v}$ will depend on $\o$ as well.  This has some subtle consequences, e.g. in Lemma \ref{lem-oreg} concerning the regularity of the value processes. The main difficulty is that we do not have a good stability result for feedback type of SDEs \reff{Xweak2}. We hope to address this issue in future research.

(iii) Note that we may get rid of the drift $b$ by using Girsanov transformation, so all  our results hold true when $b$ is random, given that  $\si$ is independent of $\o$. However, to simplify the presentation we assume $b$ is independent of $\o$ as well.
\qed}
\end{rem}

\begin{rem}
\label{rem-strong}
{\rm For each $(u, v)\in \cU_0\times \cV_0$, denote $\tilde u_t := u_t(X^{u,v})$ and $\tilde v_t := v_t(X^{u,v})$. Then $(\tilde u, \tilde v)\in \cU_0\times \cV_0$ and $\dbP^{u,v} = \dbP^{S, \tilde u, \tilde v} := \dbP^{S, \si(\tilde u, \tilde v), \wh b(\tilde u, \tilde v)}$. Thus we have $\cY^{0,{\bf 0}, u, v}_0 = \cY^{S, \tilde u, \tilde v}_0$, where,
\beaa
\left.\ba{lll}
\dis X^{u, v}_t = \int_0^t \si(s, \tilde u_s, \tilde v_s) \Big[ dB_s + b(s, \tilde u_s, \tilde v_s)ds\Big];\\ 
\dis \cY^{S, \tilde u, \tilde v}_t = \xi(X^{u,v})+ \int_t^T f(s, X^{\tilde u, \tilde v}_\cd, \cY^{S, \tilde u, \tilde v}_s, \wh\cZ^{S, \tilde u, \tilde v}_s, \tilde u_s, \tilde v_s) ds - \int_t^T \cZ^{S, \tilde u, \tilde v}_s dB_s,
\ea\right. \dbP_0\mbox{-a.s.}
\eeaa
However, we shall emphasize that the mapping from $(u, v)$ to $(\tilde u, \tilde v)$ is in pairs, and it does not induce a mapping from $u$ to $\tilde u$ (or from $v$ to $\tilde v$). Consequently the game values defined below in strong formulation are different from the $\underline Y_0$ and $\overline Y_0$ in \reff{semi-value}: 
\beaa
\underline Y^S_0 := \sup_{\tilde u \in \cU_0} \inf_{\tilde v\in \cV_0}\cY^{S, \tilde u,\tilde v}_0, && \overline Y^S_0 :=\inf_{\tilde v\in \cV_0}  \sup_{\tilde u \in \cU_0} \cY^{S, \tilde u,\tilde v}_0.
\eeaa
Indeed, in strong formulation the above game with control against control may not have the game value, namely $\underline Y^S_0 < \overline Y^S_0$, even if Isaacs condition and comparison principle for the viscosity solutions of the corresponding Bellman-Isaacs equation hold. See the counterexample Example  \ref{eg-Buckdahn} below.
\qed}
\end{rem}

\begin{rem}
\label{rem-strategy}
{\rm (i) In standard literature, see e.g. \cite{FS} and \cite{BL}, one transforms the problem into a game with strategy type of controls. That is, let $\a : \cV_t \to \cU_t$ and $\b: \cU_t \to \cV_t$ be appropriate strategies. One considers:
\beaa
  \underline Y^{"}(t,\o) := \sup_{\a} \inf_{v}\cY^{t,\o,\a(v), v}_t, &&   \overline Y^{"}(t,\o) := \inf_{\b} \sup_{u}\cY^{t,\o,u, \b(u)}_t, 
\eeaa 
This type of control problem is in fact a principal-agent problem, see e.g. Cvitanic and Zhang \cite{CZ}. 
In Markovian framework and under appropriate conditions, one can show that $  \underline Y^{"} =   \overline Y^{"}$ and is the unique 
 solution of the corresponding Bellman-Isaacs equation. However, in this formulation the two players have nonsymmetric  informations, and the lower and  upper values are defined using different information settings. In particular, it is less convenient to define saddle point in  this formulation. 
 
 (ii) Our weak formulation actually has the feature of strategy type of controls. Indeed, consider the $(\tilde u, \tilde v)$ in Remark \ref{rem-strong} again. Roughly speaking, given $u$, then $\tilde u$ is uniquely determined by $v$, which is in turn uniquely determined by $\tilde v$. Thus $u$ can be viewed as a strategy $\a$ which maps $\tilde v$ (and $B$)  to $\tilde u$. Similarly $v$ can be viewed as a strategy $\b$ which maps $\tilde u$ (and $B$) to $\tilde v$. Compared to the strategy against control, the advantage of weak formulation is that it is control against control and the two players have symmetric information.
\qed}
\end{rem}

\begin{rem}
\label{rem-HL}
{\rm When there is only drift control, namely $\si$ is independent of $(u,v)$, our formulation reduces to the work Hamadene and Lepeltier \cite{HL}.  Under Isaacs condition, by using Girsanov transformation and comparison for BSDEs, they proved $\underline Y=\overline Y$ and the existence of saddle point. We allow for both diffusion control and drift control, and we shall prove $\underline Y=\overline Y$. However, when there is diffusion control, the comparison used in \cite{HL} fails. Consequently, we are not able to follow the arguments in \cite{HL} to establish the existence of saddle point. Indeed, with the presence of diffusion control, even for stochastic optimization problem the optimal control does not seem to exist in general. We shall instead obtain some approximate saddle point in Section \ref{sect-saddle} below.  
\qed}
\end{rem}

\section{Dynamic Programming Principle}
\label{sect-DPP}
\setcounter{equation}{0}

We start with the regularity  of $\underline Y$ and $\overline Y$ in $\o$. This property is straightforward in strong formulation. Our proof here relies heavily on our assumption that $\si$ and $b$ are independent of $\o$. As pointed out in Remark \ref{rem-weak} (ii), the problem is very subtle in general case and we hope to address it in some future research. 

\begin{lem}
\label{lem-oreg}
Let Assumptions \ref{assum-sib}, \ref{assum-F},   and \ref{assum-xi} hold. Then $\underline Y$ and $\overline Y$ are uniformly continuous in $\o$ with modulus of continuity function $C\rho_0$ for some constant $C>0$. 
Consequently,   $\underline Y$ and $\overline Y$ are $\dbF$-progressively measurable.
\end{lem}
\proof 
Let $t\in [0, T],  \o, \o'\in \O$.  For any $(u, v)\in \cU_t\times \cV_t$, denote $\D \cY := \cY^{t,\o, u,  v}- \cY^{t,\o', u,  v}, \D \cZ := \cZ^{t,\o, u,  v}- \cZ^{t,\o', u,  v}$. Then, $\dbP^{t, u, v}\mbox{-a.s.}$
\beaa
\D\cY_s &=& \xi^{t,\o}(B^t_\cd) - \xi^{t,\o'}(B^t_\cd)- \int_s^T \D\cZ_r dB_r  + \int_s^T \Big[ \a_r \D \cY_r + \D \wh\cZ_r \si(r, u_r, v_r) \b_r \\
&&+[f^{t,\o}(r, B^t, \cY^{t,\o, u,  v}_r,  \wh\cZ^{t,\o, u,  v}_r, u_r, v_r) - f^{t,\o'}(r, B^t, \cY^{t,\o, u,  v}_r, \wh \cZ^{t,\o, u,  v}_r, u_r, v_r)] \Big]dr,
\eeaa
where $\a$ and $\b$ are bounded. Apply \reff{BSDEest} on the above BSDE, one obtains
\bea
\label{oregest}
|\cY^{t,\o, u,  v}_t- \cY^{t,\o', u,  v}_t|\le C\rho_0(\|\o-\o'\|_t).
\eea
Thus
\beaa
|\ul Y_t(\o) - \ul Y_t(\o')| \le \sup_{(u, v)\in \cU_t\times \cV_t} |\cY^{t,\o, u,  v}_t- \cY^{t,\o', u,  v}_t| \le C\rho_0(\|\o-\o'\|_t).
\eeaa
Similarly one can prove the estimate for $\ol Y$.
\qed

The following  Dynamical Programming Principle is important for us. 
\begin{lem}
\label{lem-DPP}
Let Assumptions \ref{assum-sib}, \ref{assum-F}, and  \ref{assum-xi} hold true. For any $0 \le s \le t \le T$ and  $\o \in \O$ we have
\beaa
\ul{Y}_s(\o) = \sup_{u \in \cU_{s}} \inf_{v \in \cV_{s}} \cY^{s,\o,u,v}_s \big[t,\ul{Y}_t^{s,\o}\big];&&
\ol{Y}_s(\o) = \inf_{v \in \cV_{s}} \sup_{u \in \cU_{s}}  \cY^{s,\o,u,v}_s \big[t,\ol{Y}_t^{s,\o}\big].
\eeaa
\end{lem}

To prove the lemma we need a technical lemma. Its proof is standard but lengthy, and is postponed  to Appendix in order not to distract our main arguments. 
\begin{lem}
\label{lem-partition}
For any $\e > 0$ and $t \in (0,T)$, there exist disjoint sets $\{ E_i, i=1,\cds, n \} \subseteq \cF_t$ such that
\beaa
 \|\o - \o'\|_t \leq \e ~\mbox{for all}~\o, \o' \in E_i, i=1,\cds, n, 
&\mbox{and}& \sup_{(u,v) \in \cU_{0}\times \cV_{0}} \dbE^{\dbP^{0,u,v}} \big( \cap_{i =1}^n  E_i^c \big) \le \e.
\eeaa
\end{lem}

{\it Proof of Lemma \ref{lem-DPP}.}
We shall prove only the Dynamic Programming Principle for $\ul{Y}$. The proof for $\ol{Y}$ is  similar. Without loss of generality, we assume  $s=0$. That is, we shall prove:
\bea
\label{DPP1}
\ul{Y}_0 &=& \sup_{u \in \cU_{0}} \inf_{v \in \cV_{0}} \cY^{0,{\bf 0}, u,v}_0 [t,\ul{Y}_t] .
\eea

{\it Step 1.} We first prove "$\ge$".  Fix arbitrary $\e > 0$ and $u\in \cU_0$. Let $\{E_i, i=1,\cds, n\}\subset \cF_t$ be given by Lemma \ref{lem-partition}, and  fix an $\o^i\in E_i$ for each $i$.  For any $\o\in E_i$, By Lemma \ref{lem-oreg} and \reff{oregest} 
we have
\bea
\label{DPPest1}
|\ul{Y}_t(\o) - \ul{Y}_t(\o^i)| \le  C\rho_0(\e)
&\mbox{ and}&
\sup_{(u,v) \in \cU_{t} \times \cV_{t}} \big| \cY^{t,\o,u,v}_t- \cY^{t,\o^i,u,v}_t \big| \leq C \rho_0(\e). 
\eea
Let $u^i \in \cU_t$ be an $\e$-optimizer of $\ul{Y}_t(\o^i)$, that is,
\bea
\label{vi}
 \inf_{v \in \cV_{t}} \cY^{t,\o^i,u^i,v}_t + \e &\ge& \ul{Y}_t(\o^i).
 \eea
Denote  $\hat{E}_n := \cap_{i=1}^n (E_i)^c$. By  Lemma \ref{lem-cUV} (i)  we define $u^\e\in \cU_0$ by:
\bea
\label{ue}
u^\e_s(\o) :=  {u}_s(\o)\1_{[0,t)}(s) + \big[  \sum_{i=1}^n u^i_s(\o^t) \mathbf{1}_{E_i} (\o) + {u}_s(\o) \mathbf{1}_{\hat{E}_n} (\o)\big]\1_{[t,T]}(s)
\eea

Now for any $v \in \cV_0$, we have
 \beaa
\cY^{0,{\bf 0}, u^\e,v}_0  &=& \cY^{0,{\bf 0},u^\e,v}_0 \big( t, \cY^{0,{\bf 0},u^\e,v}_t \big)
= \cY^{0,{\bf 0},u^\e,v}_0\Big[t, \sum_{i=1}^n \cY^{0,{\bf 0},u^\e,v}_t \1_{E_i} + \cY^{0,{\bf 0},u^\e,v}_t\1_{\hat E^n} \Big]. 
\eeaa
Since solutions of BSDEs can be constructed via Picard iteration,  one can easily check that, for any $(u,v) \in \cU_0 \times \cV_0$,
\beaa
\cY^{0,{\bf 0},u,v}_t(\o) = \cY^{t,\o,u^{t,\o},v^{t,\o}}_t,\q  \dbP^{0,u,v} \mbox{-a.e. } \o \in \O.
\eeaa
 Then it follows from \reff{vi} and Lemma \ref{lem-paste} that, for $\dbP^{0,u^\e,v}$-a.e. $\o\in E_i$,
\beaa
 \cY^{0,{\bf 0},u^\e,v}_t(\o)&=& \cY^{t,\o, (u^\e)^{t,\o},v^{t,\o}}_t= \cY^{t,\o,u^i,v^{t,\o}}_t
   \geq \inf_{v \in \cV_t} \cY^{t,\o,u^i,v}_t\\
   & \geq& \inf_{v \in \cV_t} \cY^{t,\o^i,u^i,v}_t  - C \rho_0(\e) \geq\ul{Y}_t(\o^i)- \e - C \rho_0(\e) \\
&\geq& \ul{Y}_t(\o) - \e - C \rho_0(\e).
\eeaa
Therefore, by comparison principle of BSDEs and \reff{bound} we have
 \beaa
 \cY^{0,{\bf 0},u^\e,v}_0  &\geq&  \cY^{0,{\bf 0},u^\e,v}_0\Big[ t, \sum_{i=1}^n \ul{Y}_t \1_{E_i} - (\e + C\rho_0(\e)) +\cY^{0,{\bf 0},u^\e,v}_t \1_{\hat E^n} \Big]\\
 &\ge&  \cY^{0,{\bf 0},u^\e,v}_0\Big[ t,  \ul{Y}_t  - (\e + C\rho_0(\e)) -C \1_{\hat E^n} \Big].
 \eeaa
 Recall that $\sup_{(u,v) \in \cU_{0}\times \cV_{0}} \dbP^{0,u,v}\big( \hat{E}_n \big) \le \e$. Applying \reff{BSDEestp}  we get
 \beaa
 \cY^{0,{\bf 0},u^\e,v}_0  &\ge&     \cY^{0,{\bf 0},u^\e,v}_0\big[ t, \ul{Y}_t \big]  - C(\e + \rho_0(\e))^{1\over 2} 
= \cY^{0,{\bf 0},u,v}_0\big[ t, \ul{Y}_t \big]- C(\e + \rho_0(\e))^{1\over 2}.
\eeaa
Since $v$ is arbitrary, this implies that
\beaa
\inf_{v \in \cV_0}  \cY^{0,{\bf 0},u^\e,v}_0  &\ge & \inf_{v \in \cV_0}   \cY^{0,{\bf 0},u,v}_0\big[ t, \ul{Y}_t \big]  - C(\e + \rho_0(\e))^{1\over 2}.
\eeaa
Then
\beaa
\ul{Y}_0  &\ge &   \inf_{v\in\cV_0}   \cY^{0,{\bf 0},u,v}_0\big[ t, \ul{Y}_t \big]   - C(\e + \rho_0(\e))^{1\over 2}.
\eeaa
Sending $\e\to 0$ and by the arbitrariness of $u \in\cU_0$, we obtain
\beaa
\ul{Y}_0  &\ge &\sup_{u\in \cU_0} \inf_{v\in\cV_0} \cY^{0,{\bf 0},u,v}_0\big[ t, \ul{Y}_t \big] .
\eeaa

{\it Step 2.} We now prove "$\le$". Fix $\ol u \in \cU_0$ in the form of \reff{cUV}, with $u_{ij}$ being replaced by $\ol u_{ij}$.  It suffices to prove that
\beaa
 \inf_{v \in \cV_0}   \cY_0^{t,{\bf 0},\ol{u}, v} \le \inf_{v \in \cV_0}   \cY^{0,{\bf 0},\ol{u},v}_0 [t, \ul{Y}_t] .
\eeaa
Without loss of generality, assume $t=t_{i_0}$ for some $i_0$. Notice that $\ul Y_{t_m} = \xi$, then it suffices to prove
\bea
\label{DPP2}
\inf_{v \in \cV_0}   \cY^{0,{\bf 0},\ol{u}, v}_0 [ t_{i+1}, \ul{Y}_{t_{i+1}}] \le \inf_{v \in \cV_0} \cY_0^{0,{\bf 0},\ol{u}, v} [t_i, \ul{Y}_{t_i}]  ,  &&\mbox{for all $i$}.
\eea

We now fix $i$ and recall that $\ol u_t = \sum_{j=1}^{n_i} \ol u_{ij} \1_{E^i_j}$ for $t\in [t_i, t_{i+1})$. For any $\e>0$, let $\{E_k, k=1,\cds, K\}\subset \cF_{t_i}$ be given by Lemma \ref{lem-partition}. Denote $E^i_{jk} := E^i_j \cap E_k$ and fix an $\o^{jk}\in E^i_{jk}$ for each $(j,k)$. For any $\ol v \in \cV_0$, as in Step1 we have
\beaa
 \cY^{0,{\bf 0},\ol u,\ol v}_0 [t_i,\ul{Y}_{t_i}]  &=& \cY^{0,{\bf 0},\ol u,\ol v}_0 \Big[ t_i, \sum_{j,k=1}^{n_i, K} \ul{Y}_{t_i} \1_{E^i_{jk}} +   \ul{Y}_{t_i} \1_{\cap_{k=1}^{K} E_{k}^c}\Big] \\
 &\ge& \cY^{0,{\bf 0},\ol u,\ol v}_0  \Big[ t_i, \sum_{j,k=1}^{n_i,K} \ul{Y}_{t_i} (\o^{jk})\1_{E^i_{jk}}(\o) \Big] - C(\rho_0(\e) + \e)^{1\over 2}.
 \eeaa
 By Step 1, we see that
 \beaa
 \ul{Y}_{t_i} (\o^{jk}) \ge \sup_{u\in\cU_{t_i}} \inf_{v\in \cV_{t_i}} \cY_{t_i}^{t_i,\o^{jk},u, v} \big[t_{i+1}, \ul Y^{t_i, \o^{jk}}_{t_{i+1}}\big] \ge \inf_{v \in\cV_{t_i}}  \cY_{t_i}^{t_i,\o^{jk},\ol u_{ij}, v} \big[t_{i+1}, \ul Y^{t_i, \o^{jk}}_{t_{i+1}}\big].
 \eeaa
 Here the constant $\ol u_{ij}$ denotes the constant process. Then there exists $v^{jk}\in \cV_{t_i}$ such that
 \beaa
  \ul{Y}_{t_i} (\o^{jk}) \ge   \cY_{t_i}^{t_i,\o^{jk}, \ol u_{ij},v^{jk}} \big[t_{i+1}, \ul Y^{t_i, \o^{jk}}_{t_{i+1}}\big]  - \e.
  \eeaa
  Now define 
  \beaa
  \hat v := \bar v \1_{[0, t_i)} + \Big[\sum_{j, k=1}^{n_i, K} v^{jk}(B^{t_i}) \1_{E^i_{jk}} + \bar v \1_{\cap_{k=1}^K E_k^c}\Big]\1_{[t_i, T]}.
  \eeaa
  By Lemma \ref{lem-cUV} we have $\hat v \in \cV_0$. Then, noting that $\ol u^{t_i,\o}_t = \ol u_{ij}$ for $\o\in E_{jk}^i$ and $t\in [t_i, t_{i+1})$,
  \beaa
 \cY_0^{0,{\bf 0},\ol u,\ol v} [t_i, \ul{Y}_{t_i}] &\ge& \cY_0^{0,{\bf 0},\ol u,\ol v} \Big[ \sum_{j,k=1}^{n_i,K}  \cY_{t_i}^{t_i,\o^{jk},\ol u_{ij}, v^{jk}} \big[t_{i+1}, \ul Y^{t_i, \o^{jk}}_{t_{i+1}}\big] \1_{E^i_{jk}}(\o) \Big] - C(\rho_0(\e) + \e)^{1\over 2}\\
 &=&  \cY_0^{0,{\bf 0},\ol u,\ol v} \Big[ \sum_{j,k=1}^{n_i,K}  \cY_{t_i}^{t_i,\o^{jk},\ol u^{t_i,\o}, \hat v^{t_i,\o} } \big[t_{i+1}, \ul Y^{t_i, \o^{jk}}_{t_{i+1}}\big] \1_{E^i_{jk}}(\o) \Big] - C(\rho_0(\e) + \e)^{1\over 2}\\
 &\geq& \cY_0^{0,{\bf 0},\ol u,\ol v} \Big[ \sum_{j,k=1}^{n_i,K}  \cY_{t_i}^{t_i,\o,\ol u^{t_i,\o}, \hat v^{t_i,\o} } \big[t_{i+1}, \ul Y^{t_i, \o}_{t_{i+1}}\big] \1_{E^i_{jk}}(\o) \Big] - C(\rho_0(\e) + \e)^{1\over 2}\\
 &=&  \cY_0 ^{0,{\bf 0},\ol u, \hat v} \Big[t_{i+1},  \sum_{j,k=1}^{n_i,K}  \ul Y^{t_i, \o}_{t_{i+1}}\1_{E^i_{jk}}(\o) \Big] - C(\rho_0(\e) + \e)^{1\over 2}\\
 &\ge& \cY_0^{0,{\bf 0},\ol u, \hat v} \Big[t_{i+1},  \ul Y_{t_{i+1}}\Big] - C(\rho_0(\e) + \e)^{1\over 2}\\
 &\ge& \inf_{v\in \cV_0}  \cY_0^{0,{\bf 0},\ol u, v} \Big[ t_{i+1}, \ul Y_{t_{i+1}}\Big] - C(\rho_0(\e) + \e)^{1\over 2}.
 \eeaa  
  Send $\e\to 0$, by the  arbitrariness of $\ol v \in \cV_0$ we prove \reff{DPP2}.
\qed

\begin{rem}
\label{rem-DPP}
{\rm If we use strong formulation with control against control, as in Remark \ref{rem-strong}, we can only prove the following partial Dynamic Programming Principle:
\beaa
\ul{Y}^S_s(\o) \le \sup_{u \in \cU_{s}} \inf_{v \in \cV_{s}} \cY^{s,\o,\dbP^{S,s,u,v}} \big[\ul{Y}^S_t\big];\q
\ol{Y}^S_s(\o) \ge \inf_{v \in \cV_{s}} \sup_{u \in \cU_{s}}  \cY^{s,\o,\dbP^{S,s,u,v}} \big[\ol{Y}^S_t\big],
\eeaa
which does not lead to the desired viscosity property. 
That is why we use weak formulation instead of strong formulation.
\qed}
\end{rem}

We now turn to the regularity of $\ul{Y}$ and $\ol{Y}$ in $t$, which is required for studying their viscosity property.

\begin{lem}
\label{lem-treg}
Let Assumptions \ref{assum-sib}, \ref{assum-F}, and  \ref{assum-xi} hold.  Then,  for any $0\le t_1<t_2\le T$ and  $\o\in\O$,
\bea
\label{USC1}
|\ul{Y}_{t_1}(\o) - \ul{Y}_{t_2} (\o)| +  |\ol{Y}_{t_1}(\o) - \ol{Y}_{t_2} (\o)| &\le& C  \rho_1\Big(d_\infty((t_1,\o), (t_2,\o))\Big),
\eea
where $\rho_1$ is a modulus of continuity function defined by
\bea
\label{rho1}
\rho_1(\d) := \rho_0(\d + \d^{1\over 4}) + \d+\d^{1\over 4}. 
\eea
\end{lem}
\proof We shall only prove the regularity of  $\ul{Y}$ in $t$.  The estimate for $\ol Y$ can be proved similarly.  Denote $\d := d_\infty((t_1,\o), (t_2,\o))$.

By Theorem \ref{thm-DPP} and Lemma \ref{lem-oreg} we have
\bea
\label{treg}
\Big|\ul{Y}_{t_1}(\o) - \ul{Y}_{t_2} (\o) \Big|&=& \Big| \sup_{u \in \cU_{t_1}} \inf_{v \in \cV_{t_1}}  \cY_{t_1}^{t_1,\o,u,v} \big[t_2, \ul{Y}_{t_2}^{t_1,\o}\big] - \ul{Y}_{t_2}(\o)\Big|\nonumber\\
&\leq& \sup_{u \in \cU_{t_1}, v \in \cV_{t_1}} \big |\cY_{t_1}^{t_1,\o,u,v} \big[t_2, \ul{Y}_{t_2}^{t_1,\o}\big] - \ul{Y}_{t_2}(\o) \big |.
\eea
Denote 
\beaa
\cY_t := \cY_{t}^{t_1,\o,u,v} \big[t_2, \ul{Y}_{t_2}^{t_1,\o}\big] - \ul{Y}_{t_2}(\o), && \cZ_t := \cZ_{t}^{t_1,\o,u,v} \big[t_2, \ul{Y}_{t_2}^{t_1,\o}\big].
\eeaa
Then,  $\dbP^{t_1,u,v}\mbox{-a.s.}$
\beaa
\cY_t &=& \ul{Y}_{t_2}^{t_1,\o} - \ul{Y}_{t_2}(\o) + \int_t^{t_2} f^{t_1,\o}(s, B^{t_1}, \cY_s +  \ul{Y}_{t_2}(\o),\wh \cZ_s, u_s, v_s) ds  -\int_t^{t_2} \cZ_s dB^{t_1}_s.
\eeaa
Recall from \reff{bound} that $\ul Y$ is bounded. Apply \reff{BSDEestp}  and Lemma \ref{lem-oreg}, we get
\beaa
|\cY_{t_1}| &\le &C \Big(\dbE^{\dbP^{t_1,u,v}}\big[|\ul{Y}_{t_2}^{t_1,\o} - \ul{Y}_{t_2}(\o)|^2\big]\Big)^{1\over 2} + C\d\\
&\le&C \Big(\dbE^{\dbP^{t_1,u,v}}\big[\rho^2_0(d_\infty((t_2, \o), (t_2, \o\otimes_{t_1} B^{t_1})))\big]\Big)^{1\over 2}+ C\d.
\eeaa
Note that
\beaa
&&\dbE^{\dbP^{t_1,u,v}}\big[\rho^2_0(d_\infty((t_2, \o), (t_2, \o\otimes_{t_1} B^{t_1})))\big]
\le \dbE^{\dbP^{t_1,u,v}}\Big[\rho_0^2(\d + \|B^{t_1}\|_{t_2})\big]\Big]\\
&\le& \rho^2_0(\d + \d^{1\over 4}) + C\dbP^{t_1,u,v}\Big[\|B^{t_1}\|_{t_2} \ge \d^{1\over 4}\Big] \le  \rho^2_0(\d + \d^{1\over 4}) + C\d^{-{1\over 2}}\dbE^{\dbP^{t_1,u,v}}[\|B^{t_1}\|_{t_2}^2]\\ 
&\le&  \rho^2_0(\d + \d^{1\over 4}) + C\d^{1\over 2}. 
\eeaa
Then
\beaa
|\cY_{t_1}| &\le& C\Big[\rho_0(\d + \d^{1\over 4}) + \d^{1\over 4} +\d\Big]  =C\rho_1(\d).
\eeaa 
Plug this into \reff{treg} we complete the proof.
\qed

Combining Lemmas \ref{lem-DPP} and \ref{lem-treg}, it follows from standard arguments that
\begin{thm}
\label{thm-DPP}
Let Assumptions \ref{assum-sib}, \ref{assum-F}, and  \ref{assum-xi} hold true. For any $(t,\o)\in \L$ and $\t\in \cT^t$,  we have
\beaa
\ul{Y}_t(\o) = \sup_{u \in \cU_{t}} \inf_{v \in \cV_{t}} \cY^{t,\o,u,v}_s \big[\t,\ul{Y}_\t^{t,\o}\big];&&
\ol{Y}_t(\o) = \inf_{v \in \cV_{t}} \sup_{u \in \cU_{t}}  \cY^{t,\o,u,v}_s \big[\t,\ol{Y}_\t^{t,\o}\big].
\eeaa
\end{thm}


\section{Viscosity solution properties}
\label{sect-Viscsol}
\setcounter{equation}{0}
 Define
\bea
\label{Hsemi}
\left.\ba{lll}
\dis \ul{G}(t,\o,y,z,\g) := \sup_{u \in \dbU} \inf_{v \in \dbV} \Big [ \frac{1}{2} \sigma^2(t,u,v):\g + b\si(t,u,v)z + f(t, \o, y, z\si(t, u, v), u,v)\Big]\\
\dis \ol{G}(t,\o,y,z,\g) := \inf_{v \in \dbV} \sup_{u \in \dbU} \Big [ \frac{1}{2} \sigma^2(t,u,v):\g + b\si(t,u,v)z + f(t, \o, y, z\si(t,u,v), u,v) \Big];
\ea\right.
\eea 
and consider the following path dependent PDEs:
\bea
\label{PPDE-}
-\pa_t Y_t - \ul{G}(t,\o,Y_t, \pa_\o Y_t,\pa^2_{\o\o} Y_t) = 0; \\
\label{PPDE+}
-\pa_t Y_t - \ol{G}(t,\o,Y_t, \pa_\o Y_t,\pa^2_{\o\o} Y_t) = 0.
\eea

\begin{thm}
\label{thm-viscosity}
Let Assumptions \ref{assum-sib}, \ref{assum-F}, and  \ref{assum-xi} hold. Then $\ul{Y}$ (resp. $\ol Y$) is a viscosity solution of PPDE \reff{PPDE-} (resp.  \reff{PPDE+}).
\end{thm}

\proof  We shall only prove that $\ul{Y}$ is a viscosity solution of the PPDE ($\ref{PPDE-}$). The other statement can be proved similarly. 

{\it Step 1.} We first prove the viscosity  supersolution property. Assume by contradiction that there exists $(t,\o)$ and $\f \in \ol{\cA}^L\ul{Y}(t,\o) $ such that
\beaa
c &:=& \pa_t\f(t,\mathbf{0}) + \sup_{u \in \dbU} \inf_{v \in \dbV} \big \{ \frac{1}{2} \sigma^2(t,u,v): \partial^2_{\o\o} \f(t,\mathbf{0}) + b\si(t,u,v)\partial_{\o} \f(t,\mathbf{0}) \\
&& + f(t, \o ,\ul Y_t(\o),\partial_{\o} \f(t,\mathbf{0})\si(t, u, v), u,v) \big\} > 0.
\eeaa
By Remark (\ref{rem-viscosity1}), we can assume $L$ is large enough as we will see later. Then there exists $\tilde{u} \in \dbU$ such that, for all $v \in \dbV$
\bea
\label{tildeu}
 && \pa_t\f(t,\mathbf{0}) + \frac{1}{2} \sigma^2(t,\tilde{u},v) :\partial^2_{\o\o} \f(t,\mathbf{0}) +  b\si(t,\tilde u,v)\partial_{\o} \f(t,\mathbf{0}) \notag \\
 && + f(t, \o ,\ul Y_t(\o),\partial_{\o} \f(t,\mathbf{0})\si(t, \tilde u, v), \tilde{u},v) \geq  \frac{c}{2}
\eea

Let $\ch\in \cH^t$ be the hitting time corresponding to $\f$ in \reff{cA}. For any $\e > 0$, set
\beaa
\ch_\e &:=& \inf \big \{ s \geq t : s-t+ |B^t_s| = \e \big \}.
\eeaa
By choosing $\e>0$ small enough, we have $\ch_\e\le \ch$. 
Since $\f\in C^{1,2}(\L^t)$, there exist some constant $C_\f\ge C_0$ and modulus of continuity function $\rho_\f\ge \rho_1$, which may depend on $\f$,  such that
\bea
\label{rho}
|\psi(s, B^t)| \le C_\f, ~ |\psi(s, B^t)- \psi(t, {\bf 0})| \le \rho_\f(\e),~\mbox{for} ~ t\le s \le \ch_\e, ~\psi = \f, \pa_t\f, \pa_\o \f, \pa_{\o\o}^2\f.
\eea

Now set $u := \tilde{u} \in \cU_{t}$ be a constant process and let $v \in \cV_{t}$ be arbitrary. Fix $\d>0$ and denote $\ch_{\e}^\d := \ch_\e \wedge (t+\d)$, 
\beaa
&&\cY :=  \cY^{t,\o,u,v}\big[ \ch^\d_\e, \ul{Y}^{t,\o}_{\ch_\e^\d} \big], \q  \cZ :=  \cZ^{t,\o,u,v}\big[ \ch_\e^\d, \ul{Y}^{t,\o}_{\ch_\e^\d} \big],\\
&& \D Y_s := \f(s, B^t) -\cY_s, \q \D Z_s := \pa_\o \f(s, B^t) - \cZ_s.
\eeaa
Then, applying the functional It\^{o}'s formula we obtain:
\beaa
d\D Y_s &=& \Big[\pa_t \f + {1\over 2} \pa_{\o\o}^2\f : \si^2(s, u_s, v_s) + f^{t,\o}(\cd, \cY_s, \wh\cZ_s, u_s, v_s)\Big](s, B^t) ds + \D Z_s dB^t_s\\
&=& \Big[\pa_t \f + {1\over 2} \pa_{\o\o}^2\f : \si^2(s, u_s, v_s) + f^{t,\o}(\cd, \cY_s, \wh\cZ_s, u_s, v_s) \Big](s, B^t) ds +  \D Z_s dB^t_s\\
&=& \Big[\pa_t \f + {1\over 2} \pa_{\o\o}^2\f : \si^2(s, u_s, v_s) + f^{t,\o}(\cd, \ul Y_t(\o), \pa_\o \f(\cd)\si(s, u_s, v_s), u_s, v_s)\Big](s, B^t) ds \\
&& + \Big[\a_s (\cY_s - \ul Y_t(\o)) + \D \wh Z_s \b_s\Big]ds +  \D  Z_s dB^t_s,
\eeaa
where $|\a|, |\b|\le L_0$. By \reff{tildeu} and \reff{rho}  we have
\beaa
d\D Y_s &\ge&  \Big[{c\over 2} - C_\f\rho_\f(\e)-C |\cY_s - \ul Y_t(\o)| + \D \wh Z_s \b_s\Big]ds +  \D  Z_s dB^t_s,\qq t\le s\le \ch_\e^\d.
\eeaa
Recall \reff{Wuv} and define $d\ol \dbP := M_{\ch^\d_\e} d\dbP^{t,u,v}$, where
\beaa
 M_s &:=& \exp\Big(\int_t^s [b(r, u_r, v_r) + \b_r] dW^{t,u,v}_r -{1\over 2} \int_t^s |b(r, u_r, v_r) + \b_r|^2dr\Big).
\eeaa 
Then  $\D  Z_s dB^t_s + \D \wh Z_s \b_sds$ is a $\ol \dbP$-martingale, and  thus
\beaa
\D Y_t &\le& \dbE^{\ol \dbP}\Big[\D Y_{\ch^\d_\e} - \int_t^{\ch_\e^\d}   \Big[{c\over 2} - C_\f\rho_\f(\e)-C |\cY_s - \ul Y_t(\o)|\Big]ds\Big]
\eeaa
By choosing $L$ large enough, we see that $\ol \dbP\in \cP^t_L$.  Then it follows from the definition of $\ol{\cA}^L\ul{Y}(t,\o)$ that 
\beaa
 \dbE^{\ol \dbP}[\D Y_{\ch^\d_\e}]  = \dbE^{\ol \dbP}\Big[\f(\ch^\d_\e, B^t) - \ul Y^{t,\o}_{\ch^\d_\e}\Big] \le \f(t, {\bf 0}) -  \ul Y_t(\o).
 \eeaa
Therefore, since $b$ and $\b$ are bounded,
\bea
\label{vis-est1}
&&\ul Y_t(\o) - \cY_t \le \dbE^{\ol \dbP}\Big[ \int_t^{\ch_\e^\d}   \Big[-{c\over 2} + C_\f\rho_\f(\e)+C |\cY_s - \ul Y_t(\o)|\Big]ds\Big]\nonumber\\
&\le& [-{c\over 2} + C_\f\rho_\f(\e)] \d + C_\f\d \ol \dbP(\ch_\e \le t+\d) + C\d \dbE^{\ol\dbP}\Big[ \|\cY_\cd - \ul Y_t(\o)\|_{\ch^\d_\e}\Big]\nonumber\\
&\le&  [-{c\over 2} + C_\f\rho(\e)] \d + C_\f\d \Big(\dbP^{t,u,v}(\ch_\e \le t+\d)\Big)^{1\over 2} + C\d \Big(\dbE^{\dbP^{t,u,v}}\Big[ \|\cY_\cd - \ul Y_t(\o)\|_{\ch^\d_\e}^2\Big]\Big)^{1\over 2}.
\eea
Note that,  for $\d\le {\e\over 2}$,
\bea
\label{Heest}
\dbP^{t,u,v}\Big(\ch_\e \le t+\d\Big) &\le& \dbP^{t,u,v}\Big(\d + \|B^t\|_{t+\d} \ge \e\Big) = \dbP^{t,u,v}\Big( \|B^t\|_{t+\d} \ge {\e\over 2}\Big)\nonumber\\
&\le&{C\over \e^2} \dbE^{\dbP^{t,u,v}}\big[\|B^t\|_{t+\d}^2\big] \le {C\d\over \e^2}.
\eea
Moreover, denote $\tilde\cY := \cY - \ul Y_t(\o)$. Then
\beaa
\tilde\cY_s = \ul Y^{t,\o}_{\ch^\d_\e} - \ul Y_t(\o) + \int_s^{\ch^\d_\e} f^{t,\o}(r, B^r, \tilde\cY_r+\ul Y_t(\o), \wh \cZ_r, u_r, v_r) dr -\int_s^{\ch^\d_\e} \cZ_r dB^t_r.
\eeaa
By \reff{BSDEestp} and applying Lemma \ref{lem-treg} we obtain
\bea
\label{vis-est2}
&&\dbE^{\dbP^{t,u,v}}\Big[ \|\cY_\cd - \ul Y_t(\o)\|_{\ch^\d_\e}^2\Big]  \le C\dbE^{\dbP^{t,u,v}}\Big[|\ul Y^{t,\o}_{\ch^\d_\e} - \ul Y_t(\o) |^2\Big] + C\d\nonumber\\
&\le& C\d + C\dbE^{\dbP^{t,u,v}}\Big[ \rho_1^2\big(d_\infty((t,\o), (t+\d, \o\otimes_t B^t))\big)\Big]\nonumber\\
&\le& C\d + C\dbE^{\dbP^{t,u,v}}\Big[ \rho_1^2\big(d_\infty((t,\o), (t+\d, \o)) + \|B^t\|_{t+\d})\big)\Big] \le C\rho_2(\d),
\eea
where
\bea
\label{rho2}
\rho_2(\d) &:=& \d + \sup_{(u, v)\in \cU_t\times\cV_t} \dbE^{\dbP^{t,u,v}}\Big[ \rho_1^2\big(d_\infty((t,\o), (t+\d, \o)) + \|B^t\|_{t+\d})\big)\Big].
\eea
Plug \reff{Heest} and \reff{vis-est2} into \reff{vis-est1}, we have
\beaa
\ul Y_t(\o) - \cY_t &\le& \d\Big[ -{c\over 2} + C_\f\rho_\f (\e)+ {C_\f \d^{1\over 2}\over \e}  + C\d\rho_2^{1\over 2}(\d)\Big].
\eeaa
It is clear that $\lim_{\d\to 0} \rho_2(\d) =0$. Then by first choosing $\e$ small and then choosing $\d$ small enough, we have
\beaa
\ul Y_t(\o) - \cY_t^{t,\o,u,v}[\ch^\d_\e, \ul Y^{t,\o}_{\ch^\d_\e}] &\le& -{c\over 4}\d.
\eeaa
Since $v$ is arbitrary, we get
\beaa
\ul Y_t(\o) -  \inf_{v\in \cV_t} \cY_t^{t,\o,u,v}[\ch^\d_\e, \ul Y^{t,\o}_{\ch^\d_\e}] &\le& -{c\over 4}\d,
\eeaa
which implies further that
\beaa
\ul Y_t(\o) -  \sup_{u\in\cU_t}\inf_{v\in \cV_t} \cY_t^{t,\o,u,v}[\ch^\d_\e, \ul Y^{t,\o}_{\ch^\d_\e}] &\le& -{c\over 4}\d <0.
\eeaa
This contradicts with the dynamic programming principle Theorem \ref{thm-DPP}. Therefore, $\ul Y$ is a viscosity  supersolution of PPDE \reff{PPDE-}.


{\it Step 2.} We now prove the viscosity  subsolution property. Assume by contradiction that, for some $L$ large enough,  there exists $(t,\o)$ and $\f \in \ul{\cA}^L\ul{Y}(t,\o) $ such that
\beaa
-c &:=& \pa_t\f(t,\mathbf{0}) + \sup_{u \in \dbU} \inf_{v \in \dbV} \big \{ \frac{1}{2} \sigma^2(t,u,v) :\partial^2_{\o\o} \f(t,\mathbf{0}) + b(t,u,v)\partial_{\o} \f(t,\mathbf{0}) \big\} \\
&& + f(t, \o ,\ul Y_t(\o),\partial_{\o} \f(t,\mathbf{0}), u,v) \big\} < 0.
\eeaa
Then there exists a mapping (no measurability is involved!) $\psi: \dbU \to \dbV$ such that,  for any $u \in \dbU$, 
\bea
\label{psi}
 && \pa_t\f(t,\mathbf{0}) + \frac{1}{2} \sigma^2(t,u, \psi(u) ) : \partial^2_{\o\o} \f(t,\mathbf{0}) +  b(t,u,\psi(u))\partial_{\o} \f(t,\mathbf{0})  \notag \\
 && + f(t, \o ,\ul Y_t(\o),\partial_{\o} \f(t,\mathbf{0}), u,\psi(u)) \le  -\frac{c}{2}.
\eea

For any  $u  \in \cU_{t}$, by the structure \reff{cUV}  one can easily see that $v  := \psi(u) \in \cV_{t}$. Introduce the same notations as in Step 1, and follow almost the same arguments, we obtain
\beaa
\ul Y_t(\o) - \cY_t &\ge& \d\Big[ {c\over 2} - C_\f\rho_\f (\e)- {C_\f \d^{1\over 2}\over \e}  -C\d\rho_2^{1\over 2}(\d)\Big].
\eeaa
Again,  by first choosing $\e$ small and then choosing $\d$ small enough, we have
\beaa
\ul Y_t(\o) - \cY_t^{t,\o,u,v}[\ch^\d_\e, \ul Y^{t,\o}_{\ch^\d_\e}] &\ge& {c\over 4}\d.
\eeaa
This implies
\beaa
\ul Y_t(\o) -  \inf_{v\in \cV_t} \cY_t^{t,\o,u,v}[\ch^\d_\e, \ul Y^{t,\o}_{\ch^\d_\e}] &\ge& {c\over 4}\d.
\eeaa
Since $u$ is arbitrary, then
\beaa
\ul Y_t(\o) -  \sup_{u\in\cU_t}\inf_{v\in \cV_t} \cY_t^{t,\o,u,v}[\ch^\d_\e, \ul Y^{t,\o}_{\ch^\d_\e}] &\ge& {c\over 4}\d > 0.
\eeaa
This contradicts with the dynamic programming principle Theorem \ref{thm-DPP}. Therefore, $\ul Y$ is a viscosity  subsolution of PPDE \reff{PPDE-}.
\qed

We now assume the Isaacs condition:
\bea
\label{Isaacs}
\ul{G}(t,\o,y,z,\g) = \ol{G}(t,\o,y,z,\g)=:G(t,\o,y,z,\g),
\eea
and consider the following path dependent Isaacs equation: 
\bea
\label{IsaacsPPDE}
-\pa_t Y_t- G(t,\o, Y_t, \pa_\o Y_t, \pa^2_{\o\o} Y_t) = 0.
\eea
Our main result of the paper is:
\begin{thm}
\label{thm-main}
Let Assumptions \ref{assum-sib}, \ref{assum-F}, and  \ref{assum-xi} hold. Assume further that the Isaacs condition (\ref{Isaacs}) and the uniqueness for viscosity solutions of the PPDE \reff{IsaacsPPDE}  hold. Then $\ul{Y} = \ol{Y} =: Y$ and is the unique viscosity solution of  PPDE \reff{IsaacsPPDE}.
\end{thm}
\proof Applying Theorem \ref{thm-viscosity} and  by the uniqueness of viscosity solutions, we see immediately that $\ul{Y} = \ol{Y}$ and it is the unique viscosity solution of PPDE \reff{IsaacsPPDE}.   
\qed

\begin{rem}
\label{rem-PPDEcomparison}
{\rm (i) For the comparison principle of viscosity solutions of PPDE \reff{IsaacsPPDE}, we refer to Ekren, Touzi and Zhang \cite{ETZ1}. We shall also provide a sufficient condition in Subsection \ref{sect-comparison} below.

(ii) In Markovian framework, the PPDE \reff{IsaacsPPDE} becomes a standard PDE. Note that a viscosity solution (resp. supersolution, subsolution) in the sense of Definition \ref{defn-viscosity} is a viscosity solution (resp. supersolution, subsolution) in the standard literature. Then, assuming the comparison principle for standard viscosity solution of PDEs holds true, $Y:= \ul{Y} = \ol{Y}$ and  is the unique viscosity solution of  the Bellman-Isaacs PDE  with terminal condition $Y(T,x) = \xi(x)$
\qed}
\end{rem}

 \section{Comparison principle for viscosity solutions of PPDEs}
\label{sect-comparison}
\setcounter{equation}{0}

In this section we study  the comparison principle of PPDE \reff{IsaacsPPDE}, which clearly implies the uniqueness required in Theorem \ref{thm-viscosity}. 

We first  cite a general result from \cite{ETZ3} concerning wellposedness of PPDEs, adapting to our setting. For any $(t,\o)\in \L$, denote the following deterministic function with parameter $(t,\o)$:
\bea
\label{G}
g^{t,\o}(s, y,z, \g) &:=& G(s\wedge T, \o_{\cd\wedge t}, y,z,\g).
\eea 
For any $\e>0$ and $\eta \ge 0$, we denote $T_\eta:=(1+\eta)T$, and
 \bea\label{Oet}
  \left.\ba{lll}
 O_\e:=\{x\in\dbR^d: |x|<\e\},~~ \overline O_\e :=\{x\in\dbR^d: |x|\le\e\},~~ \pa O_\e :=\{x\in\dbR^d: |x|=\e\};\ms\\
 \cO^{\e,\eta}_t :=  [t,T_\eta)\times O_\e,~
 \overline \cO^{\e,\eta}_t 
  := [t, T_\eta]\times \overline O_\e),~
 \pa \cO^{\e,\eta}_t 
 :=\big([t,T_\eta]\times\partial O_\e\big) 
   \cup \big(\{T_\eta\}\times O_\e\big),
   \ea\right.
 \eea
Consider the following  localized and path-frozen PDE defined for every $(t,\o)\in\L$:
 \bea\label{PDEe}
 \mbox{(E)}^{t,\o}_{\eps,\eta}
 \quad
 &&
 \mathbf{L}^{t,\o}th
 :=
 -\pa_t \th - g^{t,\o} (s, \th, D  \th, D^2  \th) 
 = 0
 ~~\mbox{on}~~
 \cO^{\e,\eta}_t.
 \eea
Here $\pa_t, D, D^2$ are standard differential operators. 
\begin{assum}\label{assum-comparison}
For any $\e>0,\eta \ge 0$, $(t,\o)\in\L$, and any  $h\in C^0(\pa \cO^{\e,\eta}_t)$,
we have  $\overline \th = \underline \th$, where
 \bea\label{barv}
 \left.\ba{lll}
 \overline \th(s,x) 
 &:=& 
 \inf\Big\{w(s,x): w~\mbox{classical supersolution of 
                          {\rm (E)}}^{t,\o}_{\eps,\eta}
                     ~\mbox{and}
                     ~w \ge h 
                     ~\mbox{on}~\pa \cO^{\e,\eta}_t\Big\},\ms
 \\
 \underline \th(s,x) 
 &:=& 
 \sup\Big\{w(s,x): w~\mbox{classical subsolution of 
                          {\rm (E)}}^{t,\o}_{\eps,\eta}
                    ~\mbox{and}
                    ~w \le h 
                    ~\mbox{on}~\pa \cO^{\e,\eta}_t\Big\}.
 \ea\right.
 \eea 
\end{assum}
By  \cite{ETZ3} Theorem 3.4, we have

\begin{thm}
\label{thm-wellposed}
Let Assumptions \ref{assum-sib}, \ref{assum-F},  \ref{assum-xi}, and  the Isaacs condition (\ref{Isaacs}) hold. Then, under the additional Assumption \ref{assum-comparison}, the  PPDE \reff{IsaacsPPDE} has a unique viscosity solution and  the comparison principle of viscosity solutions holds. 
\end{thm}

We remark that Assumption \ref{assum-comparison} is in the spirit of Perron's approach. However, in standard literature the $w$ in \reff{barv} is required only to be viscosity supersolution or subsolution, while we require it to be a classical one.  To check that,  we present a result concerning classical solutions of parabolic PDEs.   

We first simplify the notations. Let $O \subset \dbR^d$ be open, connected, bounded, and with smooth boundary. Set 
\beaa
\cO := [0, T)\times O,\q  \ol \cO := [0, T]\times  \ol O,\q \pa \cO := \big( [0, T] \times \pa O\big) \cup \big( \{T\} \times O\big).
\eeaa
Consider the following (standard) PDE in $\cO$ with boundary condition $h$: 
\bea
\label{PDEg}
-\pa_t \th - g (t, x, \th, D \th, D^2\th) = 0 ~\mbox{in}~ \cO &\mbox{and}& \th = h ~\mbox{on}~\pa \cO. 
\eea
 Then we have the following result, whose argument is standard in the literature and is communicated to us by Lihe Wang. We present its proof in Appendix for completeness.

\begin{lem}
\label{lem-Wang}
Assume

(i) $h\in C^{1,2}(\ol \cO)$ and $g(\cd, y,z,\g) \in C^{1,2}(\ol\cO)$ for any $(y,z,\g)$;

(ii)  $g$ is continuously differentiable in $(y, z, \g)$ with bounded derivatives; 

(iii) $\pa_\g g \ge c_0 I_d$ for some $c_0>0$, and  $d \le 2$.

\no Then the PDE \reff{PDEg} has a classical solution $\th \in C^{1,2}(\ol \cO)$.
\end{lem}

We now have
\begin{prop}
\label{prop-Perron}
Let Assumptions \ref{assum-sib}, \ref{assum-F}, \ref{assum-xi}, and  the Isaacs condition (\ref{Isaacs}) hold.  Assume further that
\bea
\label{Wang}
\mbox{ $\si \ge c_0 I_d$ for some $c_0 >0$ and the dimension $d \le 2$.}
\eea
 Then Assumption \ref{assum-comparison} holds true.  Consequently $\ul{Y} = \ol{Y} =: Y$ and is the unique viscosity solution of  PPDE \reff{IsaacsPPDE}.
\end{prop}
\proof  We use the notations in Assumption \ref{assum-comparison}. By \cite{ETZ2} Proposition 3.14, we may assume without loss of generality that
\bea
\label{Gmonotone}
G(\cd, y_1, \cd) - G(\cd, y_2, \cd) \le y_2-y_1&\mbox{for any}& y_1 \ge y_2
\eea

First, one can easily extend $h$ to a uniformly continuous function on $[t, \infty) \times \dbR^d$, still denoted  as $h$. For any $\d>0$, let $g^{t,\o}_\d$ and $h_\d$ be smooth mollifiers of $g^{t,\o}$ and $h$ such that $\|g^{t,\o}_\d - g\|_\infty \le \d, \|h_\d - h\|_\infty \le \d$.  By our assumptions, it is clear that $c_0 I_d \le \pa_\g g^{t,\o}_\d \le L_0 I_d$. Apply Lemma \ref{lem-Wang},  the following PDE has a classical solution $\th_\d\in C^{1,2}(\ol \cO^{\e,\eta}_t)$: 
 \beaa
 -\pa_t  \th_\d - g^{t,\o}_\d(s,  \th_\d, D  \th_\d, D^2  \th_\d) =0,~~\mbox{in}~\cO^{\e,\eta}_t,\q \th_\d =  h_\d ~~\mbox{on} ~\pa \cO^{\e,\eta}_t.
  \eeaa
Denote
\beaa
\ol \th_\d := \th_\d + \d &\mbox{and}& \ul \th_\d := \th_\d - \d.
\eeaa
Then clearly $\ol\th_\d \in C^{1,2}(\ol \cO^{\e,\eta}_t)$, $\ol\th_\d \ge h$ on $\pa  \cO^{\e,\eta}_t$. Moreover, by \reff{Gmonotone}
\beaa
 \mathbf{L}^{t,\o}  \ol\th_\d
 &=& -\pa_t \th_\d - g^{t,\o} (s, \th_\d + \d, D  \th_\d, D^2  \th_\d) \\
 &\ge& -\pa_t \th_\d - g^{t,\o} (s, \th_\d, D  \th_\d, D^2  \th_\d)  + \d\\
 &=& g^{t,\o}_\d(s,  \th_\d, D  \th_\d, D^2  \th_\d) - g^{t,\o} (s, \th_\d, D  \th_\d, D^2  \th_\d)  + \d \ge 0.
\eeaa
Then $\ol\th_\d$ is a classical supersolution of ${\rm (E)}^{t,\o}_{\eps,\eta}$, and thus $\ol \th \le \ol\th_\d$. Similarly,  $\ul \th \le \ul\th_\d$. Then
\beaa
0\le \ol \th - \ul \th \le  \ol \th_\d -  \ul \th_\d = 2\d.
\eeaa
Since $\d>0$ is arbitrary, we conclude that $ \ol \th = \ul \th$.
\qed

\section{Approximate saddle point}
\label{sect-saddle} 
\setcounter{equation}{0}

In this section we discuss briefly saddle points of the game, assuming the game value exists.
In our setting, it is natural to define

\begin{defn}
\label{defn-saddle}
We call $(u^*, v^*) \in \cU_0 \times \cV_0$ a saddle point of the game if 
\beaa
\cY^{0, {\bf 0}, u, v^*}_0 \le \cY^{0, {\bf 0}, u^*, v^*}_0  \le \cY^{0, {\bf 0}, u^*, v}_0 &\mbox{for all}& u\in \cU_0, v\in \cV_0.
\eeaa
\end{defn}
We remark that, if a saddle point $(u^*, v^*)$ exists, then it is straightforward to check that the game has a value $Y_0 :=  \cY^{0, {\bf 0}, u^*, v^*}_0$. However, even in stochastic optimization problem with diffusion control, in general the optimal control may not exist. We thus study approximate saddle points only.
\begin{defn}
\label{defn-approxsaddle}
For any $\e>0$, we call $(u^\e, v^\e) \in \cU_0 \times \cV_0$ an $\e$-saddle point of the game if 
\beaa
\cY^{0, {\bf 0}, u, v^\e}_0 - \e \le \cY^{0, {\bf 0}, u^\e, v^\e}_0  \le \cY^{0, {\bf 0}, u^\e, v}_0 + \e &\mbox{for all}& u\in \cU_0, v\in \cV_0.
\eeaa
\end{defn}
We have the following simple observation:

\begin{prop}
\label{prop-saddle}
Assume the game has a value,  then it  has an $\e$-saddle point   $(u^\e, v^\e)$  for any $\e>0$. 
\end{prop}
\proof Let $\underline Y_0 = Y_0 = \overline Y_0$ be the game value. Then for any $\e>0$, there exist $u^\e\in \cU_0, v^\e \in \cV_0$ such that
\beaa
Y_0 - \e < \inf_{v\in \cV_0} \cY^{0, {\bf 0}, u^\e, v}_0  \le Y_0 \le \sup_{u\in \cU_0} \cY^{0, {\bf 0}, u, v^\e}_0 \le Y_0+\e.
\eeaa
In particular, this implies that
\beaa
Y_0 - \e < \inf_{v\in \cV_0} \cY^{0, {\bf 0}, u^\e, v}_0  \le \cY^{0, {\bf 0}, u^\e, v^\e}_0 \le \sup_{u\in \cU_0} \cY^{0, {\bf 0}, u, v^\e}_0 \le Y_0+\e.
\eeaa
That is, $(u^\e, v^\e)$ is an $\e$-saddle point. Moreover, we observe that $|\cY^{0, {\bf 0}, u^\e, v^\e}_0-Y_0|\le \e$.
\qed


\section{Appendix}
\label{sect-Appendix}
\setcounter{equation}{0}

\subsection{ Proof of Lemma \ref{lem-paste}}
We prove the case $\overline\cA^S$ only. The case $\overline\cA^W$ can be proved similarly. Let $X$ be the unique strong solution to SDE \reff{weaksolution}  with coefficients $(\si, b)$, and $X^i$ be the unique strong solution to SDE \reff{weaksolution} on $[t, T]$ with coefficients $(\si^i, b^i)$. 

First, denote 
\beaa
\bar X_s &=& X_s \1_{[0, t)}(s) + \Big[X_t + \sum_{i=1}^n \1_{E_i}(X) X^i_s(B^t)\Big] \1_{[t, T]}(s),\q 0\le s\le T.
\eeaa 
One can check straightforwardly that $\bar X$ is a strong solution to SDE \reff{weaksolution} with coefficients $(\bar \si, \bar b)$. On the other hand, let $\tilde X$ be an arbitrary  strong solution to SDE \reff{weaksolution} with coefficients $(\bar \si, \bar b)$. Then both $\bar X$ and $\tilde X$ satisfy SDE  \reff{weaksolution} on $[0, t]$ with coefficients $(\si, b)$. By the uniqueness assumption of $(\si, b)$, we see that $\bar X = \tilde X$ on $[0, t]$, $\dbP_0$-a.s. In particular, this implies $\1_{E_i}(\bar X) = \1_{E_i}(\tilde X)$.  Then for $\dbP_0$-a.e. $\o\in \O$, there exists unique $i$ such that $\1_{E_i}(\bar X) = \1_{E_i}(\tilde X)=1$. Thus both $\bar X^{t,\o}$ and $\tilde X^{t,\o}$ satisfy SDE \reff{weaksolution} on $[t, T]$ with coefficients $(\si^i, b^i)$. By the uniqueness assumption of $(\si^i, b^i)$, we see that $\bar X^{t,\o}=\tilde X^{t,\o}$, $\dbP^t_0$-a.s. This implies that $\bar X = \tilde X$, $\dbP_0$-a.s. and therefore, $(\bar \si, \bar b)\in \ol \cA^S$.

Finally, since $\bar X = X$ on $[0, t]$, we have $\dbP^{\bar\si, \bar b} = \dbP^{\si, b}$ on $\cF_t$. Moreover, since $\bar X^{t,\o}(B^t) = X_t(\o)+X^i(B^t)$ whenever $\1_{E_i}(X) =1$, by the definition of r.c.p.d. we see that $(\dbP^{\bar\si, \bar b} )^{t, \o}= \dbP^{t,\si^i, b^i}$ for $\dbP^{\si, b}$-a.e. $\o\in E_i$.
\qed

\subsection{Proof of Lemma \ref{lem-BSDE}}

Recall the $\dbP^{t,u,v}$-Brownian motion $W^{t,u,v}$ defined in \reff{Wuv}. One may rewrite BSDE \reff{BSDEuv} as
\beaa
\cY_s = \eta + \int_s^{\t} \Big[f^{t,\o}(r,B^{t}_\cd, \cY_r, \wh\cZ_r, u_r, v_r) + \wh \cZ_r b(r, u_r, v_r)\Big]dr - \int_s^{\t} \wh\cZ_r dW^{t,u,v}_r, \; \dbP^{t,u,v} \as
\eeaa
Then \reff{BSDEest} follows from standard BSDE arguments. Moreover, note that
\beaa
\cY_s = \eta + \int_s^{\t} \Big[f^{t,\o}(r,B^{t}_\cd, 0, {\bf 0}, u_r, v_r) + \a_r \cY_r + \wh \cZ_r \b_r\Big]dr - \int_s^{\t} \wh\cZ_r dW^{t,u,v}_r, \; \dbP^{t,u,v} \as
\eeaa
where $\a, \b$ are bounded. Denote
\beaa
\G_r &:=& \exp\Big(\int_t^r \b_s dW^{t,u,v}_s  + \int_t^r [\a_r - {1\over 2} |\b_r|^2]dr\Big).
\eeaa
Then
\beaa
\cY_t = \G_\t \eta + \int_s^{\t} \G_r f^{t,\o}(r,B^{t}_\cd, 0, {\bf 0}, u_r, v_r) dr - \int_s^{\t} [\cds]  dW^{t,u,v}_r, \; \dbP^{t,u,v} \as
\eeaa
 Thus
\beaa
|\cY_t| &=& \Big|\dbE^{\dbP^{t,u,v}}\Big[ \G_\t \eta + \int_s^{\t} \G_r f^{t,\o}(r,B^{t}_\cd, 0, {\bf 0}, u_r, v_r) dr\Big]\Big| \\
&\le& \Big(\dbE^{\dbP^{t,u,v}}[ \G_\t^2]\Big)^{1\over 2}\Big(\dbE^{\dbP^{t,u,v}}[ |\eta|^2]\Big)^{1\over 2}\\
&&+ \d  \Big(\dbE^{\dbP^{t,u,v}}[\|\G\|_\t^2]\Big)^{1\over 2}\Big(\dbE^{\dbP^{t,u,v}}\big[ \int_s^{\t} |f^{t,\o}(r,B^{t}_\cd, 0, {\bf 0}, u_r, v_r)|^2 dr \big]\Big)^{1\over 2}.
\eeaa
It is clear that $\dbE^{\dbP^{t,u,v}}[\|\G\|_\t^2]\le C$. Then \reff{BSDEestp} follows immediately.
\qed

\subsection{ Proof of Lemma \ref{lem-partition}}
We introduce the following capacity $\cC$:
\bea
\label{cC}
\cC(A) := \sup_{(u, v)\in \cU_0\times \cV_0} \dbP^{0,u,v}(A), &\mbox{for all}& A\in \cF_T.
\eea
In this proof we abuse a notation a little bit by denoting $B^s_r := B_r - B_s$ for $0\le s\le r\le t$.

{\it Step 1.} We first show that,  for any $c, \d > 0$, and $R>0$,
\bea
\label{cCest1}
\cC\Big(\|B\|_t > R\Big) \le {C\over R^4} &\mbox{and}& \cC\Big(\sup_{0\le s\le t} \|B^s\|_{(s+\d)\wedge t}\ge c\Big) \le {C\d\over c^4}.
\eea
Indeed,  for any $(u, v)\in \cU_0\times\cV_0$ and any $0\le t_1<t_2$, since $\si$ and $b$ are bounded, then by \reff{Wuv} and applying the Burkholder-Davis-Gundy Inequality we get
\bea
\label{cCest2}
\dbE^{\dbP^{0,u,v}}\Big[\|B^{t_1}\|_{t_2}^4\Big]
& =& \dbE^{\dbP^{0,u,v}}\Big[\sup_{t_1\le s\le t_2}\Big|\int_{t_1}^{t_2} \sigma(r,u_r,v_r) b(r, u_r, v_r) dr + \int_{t_1}^{t_2} \sigma(r,u_r,v_r) dW^{u,v})_r\Big|^4\Big]\nonumber\\
&\le& C\dbE^{\dbP^{0,u,v}}\Big[\Big(\int_{t_1}^{t_2}|\sigma(r,u_r,v_r) b(r, u_r, v_r)| dr\Big)^4 + \Big(\int_{t_1}^{t_2}|\sigma(r,u_r,v_r)|^2 d_r\Big)^2\Big]\nonumber\\
&\le& C(t_2-t_1)^2.
\eea
Then
\beaa
\dbP^{0,u,v}\Big(\|B\|_t > R\Big)\le {1\over R^4}\dbE^{\dbP^{0,u,v}}\Big[\|B\|_t^4\Big] \le {C\over R^4}.
\eeaa
By the definiton of $\cC$, this implies the first estimate in \reff{cCest1}.

Next, let $0=t_1<\cds<t_m = t$ such that $\d \le \D t_i < 2\d$ for all $i$.  Then
\beaa
\sup_{0\le s\le t} \|B^s\|_{(s+\d)\wedge t}&=& \max_{0\le i\le m-1} \sup_{t_i\le s\le t_{i+1}} \sup_{s\le r \le (s+\d)\wedge t}|B_r-B_s| \\
&\le&   \max_{0\le i\le m-1} \sup_{t_i\le s\le t_{i+1}} \sup_{s\le r \le (s+\d)\wedge t}\Big[|B_r-B_{t_i}| + |B_s-B_{t_i}|\Big] \\
&\le& 2  \max_{0\le i\le m-1} \|B^{t_i}\|_{t_{i} + 3\d} .\eeaa
Then, noting that $m\le {T\over\d}$,  by \reff{cCest2} we have
\beaa
\dbP^{0,u,v}\Big(\sup_{0\le s\le t} \|B^s\|_{(s+\d)\wedge t}\ge c\Big) &\le& {1\over c^4}\dbE^{\dbP^{0,u,v}}\Big[\sup_{0\le s\le t} \|B^s\|_{(s+\d)\wedge t}^4 \Big] \\
&\le& {C\over c^4}\sum_{i=0}^{m-1}\dbE^{\dbP^{0,u,v}}\Big[\|B^{t_i}|^4_{t_i+3\d} \Big] \le {C\over c^4} m \d^2 \le {C\d\over c^4}.
\eeaa
By the definition of $\cC$ we obtain the second estimate in  \reff{cCest1}.

{\it Step 2.} We now fix $\e>0$. For the constant $C$ in \reff{cCest1}, set 
\beaa
c := {\e\over 3}, \q \d :=  {c^4\e \over 2C}\wedge t = {\e^5\over 162 C}\wedge t,\q R:= ({2C \over \e})^{1\over 4}.
\eeaa
 Let $0=t_0<\cds<t_m = t$ such that $\d\le \D t_i \le 2\d$, $i=1,\cds, m$. Clearly there exists a partition $\{\tilde E_j, 1\le j\le n\} \subset \cF_t$ such that
 \beaa
 \cup_{j=1}^n \tilde E_j = \Big\{ \max_{0\le i\le m} |B_{t_i}| \le R+c\Big\} &\mbox{and}& \max_{0\le i\le m} |\o_{t_i} - \o'_{t_i}| \le {\e\over 3} ~\mbox{for all}~\o, \o'\in \tilde E_i.
 \eeaa
 Now set 
 \beaa
 E_j:= \tilde E_j \cap A, &\mbox{where}& A:=\Big\{\sup_{0\le s\le t} \|B^s\|_{(s+\d)\wedge t}\le c\Big\}\in\cF_t. 
 \eeaa
 Then for any $\o, \o' \in E_j$,
 \beaa
 \|\o-\o'\|_t &=&\max_{0\le i \le m-1} \sup_{t_i \le s\le t_{i+1}}|\o_s - \o'_s| \\
 &\le& \max_{0\le i \le m-1} \sup_{t_i \le s\le t_{i+1}}\Big[|\o_s - \o_{t_i}| + |\o'_s- \o'_{t_i}| + |\o_{t_i} - \o'_{t_i}|\Big]\\
 &\le& \max_{0\le i \le m-1} \sup_{t_i \le s\le t_{i+1}}\Big[{\e\over 3} +{\e\over 3} +{\e\over 3}  \Big] = \e.
 \eeaa
 On the other hand, 
 \beaa
 \cap_{j=1}^n E_j^c &=& \Big(\cup_{j=1}^n \tilde E_j\Big)^c \cup A^c =  \Big\{ \max_{0\le i\le m} |B_{t_i}| > R+c\Big\} \cup  A^c\\
 & \subset& \Big(\Big\{ \max_{0\le i\le m} |B_{t_i}| > R+c\Big\} \cap A\Big) \cup A^c.
  \eeaa
  For each $\o\in\Big\{ \max_{0\le i\le m} |B_{t_i}| > R+c\Big\} \cap A$, we have
 \beaa
 \|\o\|_t = \max_{0\le i \le m-1} \sup_{t_i \le s\le t_{i+1}}|\o_s| \ge \max_{0\le i \le m-1} \sup_{t_i \le s\le t_{i+1}}\Big[|\o_{t_i}| - |\o_s- \o_{t_i}|\Big] > (R+c)-c =R.
 \eeaa
 That is,
 \beaa
 \Big\{ \max_{0\le i\le m} |B_{t_i}| > R+c\Big\} \cap A &\subset& \Big\{\|B\|_t > R\Big\},
 \eeaa
 and therefore,
 \beaa
  \cap_{j=1}^n E_j^c \subset \Big\{\|B\|_t > R\Big\} \cup A^c.
  \eeaa
  Now it follows from \reff{cCest1} that $\cC(\cap_{j=1}^n E_j^c) \le \e$.
  \qed

\subsection{Proof of Lemma \ref{lem-Wang}}
\label{sect-Wang}
As standard in PDE literature, it suffices to provide a priori estimates. That is, we assume $\th\in C^{1,2}(\ol\cO)$ satisfies PDE \reff{PDEg}, and we shall provide estimates which depends only on the parameters in our assumptions. 

(i) We first establish the estimates in the case $g= g(\g)$. We proceed in several steps.
 
 \ms
{\it Step 1.} We first cite a result from Ladyzenskaya et al  \cite{ LS}.  Assume $\th$ satisfying the following linear PDE:
\beaa
-\pa_t \th  - {1\over 2} A(t,x) : D^2  \th =0,
\eeaa
 where $ A = [a_{ij}]_{1\le i, j \le d}$ is required only to be measurable and ${\bf 0} < c_0 I_d \le A \le C_0 I_d$ . Then $\th \in C^{{\a\over 2},\a}_{\rm loc}$, where $\a$ depends only on $c_0$ and $C_0$. 
 
 \ms
{\it Step 2.}  We next cite a result by  Caffarelli \cite{Caffarelli}.  
\beaa
&\mbox{The elliptic PDE $g(D^2 \th) = f(x)$ with $f\in C^\a$ has $C^{2,\a}$-solution}&\\
&\mbox{ if the simplified PDE $g(D^2\th) = $ constant has $C^{2,\tilde \a}$-solution for some $\tilde \a>\a$.}&
\eeaa 
 
\ms
{\it Step 3.}   We also need the DeGiorgi-Nash estimate: If $\sum_{i,j=1}^d D_{x_i} (a_{ij} D_{x_j} \th) = 0$, then $\th\in C^\a$. See, e.g., Gilbarg and Trudinger \cite{GT} Theorem 8.22.

\ms

{\it Step 4.} We now come back to  the PDE \reff{PDEg} with $g=g(\g)$.  First, set $\tilde \th := \pa_t \th$. Differentiate both sides of \reff{PDEg} with respect to $t$ we obtain:
\beaa
\pa_t \tilde \th +  [\pa_{\g} g (D^2 \th)] : D^2 \tilde \th = 0.
\eeaa
By Step 1, we have $\pa_t \th = \tilde \th \in C^{{\a\over 2},\a}$.  Now  fix $t$. Then \reff{PDEg} becomes
\beaa 
 g(D^2 \th) = - \pa_t \th \in C^\a.
 \eeaa
 By Step 2, it suffices to show that
 \bea
 \label{PDEconstant}
 g(D^2 \th) = \mbox{constant  has $C^{2,\a'}$-solution for some $\a'>\a$}
 \eea
For this, we can only prove in the cases $d=1$ or $d=2$.
 
 In the case $d=1$, notice that $g$ is strictly increasing, then $D^2 \th =$ constant and thus $\th$ is a parabola. 
 
 In the case $d=2$, fix $k=1,2$ and denote $\th^k := D_{x_k} \th$. Differentiate both sides of \reff{PDEconstant} with respect to $x_k$:
 \beaa
A: D^2 \th^k =0, &\mbox{where}& A := [a_{i, j}]_{1\le i, j\le 2} :=\pa_{\g} g(D^2 \th)
 \eeaa
 Note that $a_{11} \ge c_0>0$. Then
 \beaa
 D^2_{x_1x_1} \th^k + {a_{12}\over a_{11}} D^2_{x_1x_2} \th^k + {a_{22}\over a_{11}} D^2_{x_2x_2} \th^k =0.
 \eeaa
 For $l=1,2$, differentiate both sides of the above PDE with respect to $x_l$ and denote $\th^{k,l} := D_{x_l} \th^k = D^2_{x_kx_l} \th$:
 \beaa
 D^2_{x_1x_1} \th^{k,l} + D_{x_l}\Big({a_{12}\over a_{11}} D^2_{x_1x_2} \th^k\Big) +  D_{x_l}\Big({a_{22}\over a_{11}} D^2_{x_2x_2} \th^k\Big) =0.
 \eeaa
 In the case $l=2$, this is:
 \beaa
 D_{x_1} (D_{x_1} \th^{k,2}) + D_{x_2}\Big({a_{12}\over a_{11}} D_{x_1} \th^{k,2}\Big) +  D_{x_2}\Big({a_{22}\over a_{11}} D_{x_2} \th^{k,2}\Big) =0.
 \eeaa
 By Step 3, $\th^{k,2} \in C^\a$. Similarly, $\th^{k,1}\in C^\a$. That is, for any $t$, $\th(t,\cd) \in C^{2+\a}$. Moreover, it follows from PDE \reff{PDEg} that $\th$ is differentiable in $t$ and thus $\th \in C^{1,2}$.

{\it (ii).} We now consider the general case where $g = g (t,x, y, z, \g)$. We define a map $J: C^{1,2}(\ol\cO) \rightarrow C^{1,2}(\ol\cO)$ by $J\th := \tilde\th$, where, thanks to (i),  $\tilde\th$ is the classical solution of the following PDE:
\beaa
 -\pa_t \tilde\th - g (t, x, \th, D \th, D^2 \tilde\th) 
 = 0~\mbox{in}~\cO &\mbox{and}&  \tilde\th = \th ~\mbox{on}~\pa \cO.
\eeaa
Now, following the arguments in \cite{Lieberman} Theorem 8.2, one can show that the mapping $J$ is a contraction mapping if $T$ is small enough. Moreover, the fixed point $\th$ of the mapping $J$ is also in $C^{1,2}(\ol\cO)$. Therefore, we can conclude the so called small time existence: the PDE \reff{PDEg} has a classical solution when $T$ is small enough. 

Next, \cite{Lieberman} Theorem 14.4 gives an a priori uniform estimate for the H\"{o}lder-$(1+\d)$ norm of the classical solution to \reff{PDEg}, for some $\d \in (0,1)$, where the definition of  the H\"{o}lder-$(1+\d)$ norm is given in \cite{Lieberman} Chapter IV, Section 1. Using this a priori estimate and following the arguments  in \cite{Lieberman} Theorem 8.3, we can infer the existence of the classical solution over arbitrary time duration $[0, T]$  from the small time existence, and thus complete the proof.
\qed

\subsection{Buckdahn's counterexample}

As pointed out in Remark \ref{rem-strong}, a game with control against control in strong formulation may not have the game value, even if the Isaacs condition and the comparison principle for the associate Bellman-Isaacs equation hold.  The following  counterexample is communicated to us by Rainer Buckdahn.

\begin{eg}
\label{eg-Buckdahn}
Let $d=2$, $\dbU := \{x\in \dbR: |x|\le 1\}$, $\dbV := \{x\in \dbR: |x|\le 2\}$, and  $\cU$ (resp. $\cV$) be the set of $\dbF$-progressively measurable $\dbU$-valued (resp. $\dbV$-valued) processes. Write $B = (B^1, B^2)$. Given $(u, v) \in \cU\times \cV$, the controlled state process $X^{u,v} = (X^{1,u}, X^{2, v})$ is determined by:
\beaa
X^{1, u}_t := \a B^1_t + \int_0^t u_s ds, && X^{2,  v}_t := \a B^2_t + \int_0^t v_s ds
\eeaa
where $\a\ge 0$ is a constant.  Define, for some $a\in \dbR$,
\beaa
J(u,v) := \dbE^{\dbP_0}\Big[|a + X^{1, u}_T-X^{2, v}_T|\Big],\q  \ul Y_0 := \sup_{u\in \cU}\inf_{v\in \cV}  J(u,v),\q \ol Y_0 :=\inf_{v\in \cV} \sup_{u\in \cU} J(u,v).
\eeaa
Then, for  $0\le \a < \sqrt{T\over 2}$ and $|a|\le T$, we have $\ul Y_0 < \ol Y_0$.
\end{eg}

\proof For any $ u \in \cU$, set $ v_t:=  u_t + {a \over T}$. Then $ v \in \cV$ and,
\beaa
a+ X^{1, u}_T - X^{2, v}_T = a + \a B^1_T + \int_0^T  u_t dt - \a B^2_T - \int_0^T [ u_t + {a\over T}]dt = \a[B^1_T - B^2_T].
\eeaa
Thus
\beaa
J(u, v) = \a \dbE^{\dbP_0}\Big[|B^1_T-B^2_T|\Big] = \a \sqrt{2T}. 
\eeaa
This implies that 
$
\inf_{v\in \cV} J( u,  v)   \le \a \sqrt{2T}.
$
Since $u$ is arbitrary, we get
\bea
\label{Buckdahn1}
\ul Y_0 \le \a \sqrt{2T}.
\eea

On the other hand, for any $ v \in \cV$, set
\bea
\label{Buckdahn2}
u_t :=  u_0 :=  {a-\dbE^{\dbP^0}[X^{2,v}_T]\over |a-\dbE^{\dbP^0}[X^{2,v}_T]|} \1_{\{a-\dbE^{\dbP^0}[X^{2,v}_T]\neq 0\}} + \1_{\{a-\dbE^{\dbP^0}[X^{2,v}_T] = 0\}}.
\eea
That is, $ u$ is a constant process. One can easily check that 
\beaa
 u\in \cU,\q  |u_0|=1, \q a-\dbE^{\dbP^0}[X^{2,v}_T] = u_0 |a-\dbE^{\dbP^0}[X^{2,v}_T]|.
\eeaa
Then
\beaa
\dbE^{\dbP_0}\Big[a+ X^{1, u}_T - X^{2,  v}_T \Big]= a+  u_0 T  - \dbE^{\dbP_0}[  X^{2, v}_T] =  u_0\Big[T+  |a-\dbE^{\dbP^0}[X^{2,v}_T]|\Big].
\eeaa
Thus,
\beaa
J(u,v) &\ge& \Big|\dbE^{\dbP_0}\Big[a+ X^{1, u}_T - X^{2, v}_T\Big]\Big|= |u_0|\Big[T+  |a-\dbE^{\dbP^0}[X^{2,v}_T]|\Big]\\
&=& T+  |a-\dbE^{\dbP^0}[X^{2,v}_T]| \ge T.
\eeaa
This implies $\sup_{u\in\cU} J(u, v) \ge T$. Since $v$ is arbitrary, we have $\ol Y_0 \ge T$. This, together with \ref{Buckdahn1}, implies that $\ul Y_0 < \ol Y_0$ when $0\le \a< \sqrt{T\over 2}$.  

Moreover,   note that in this case the system is Markovian and $\pa_\o Y =D Y$. The Hamiltonians in \reff{Hsemi} become: for $(t,x,y,z,\g) \in [0,T]\times \dbR^2 \times \dbR \times \dbR^2 \times \dbS^2$,
\beaa
\ul G(t,x,y,z, \g) := \sup_{u\in \dbU} \inf_{v\in\dbV}\Big[{1\over 2} \a\tr(\g) + u z_1 + v z_2\Big] = {1\over 2} \a \tr(\g) + z_1^+ - 2z_2^-;\\
 \ol G(t,x,y,z, \g) :=  \inf_{v\in\dbV}\sup_{u\in \dbU}\Big[{1\over 2} \a\tr(\g) + u z_1 + v z_2\Big] = {1\over 2} \a \tr(\g) + z_1^+ - 2z_2^-.
\eeaa
Then the Isaacs condition holds, and the corresponding Bellman-Isaacs equation becomes:
\beaa
- \pa_t Y_t - {1\over 2} \a  \Big[D^2_{x_1x_1} Y_t +D^2_{x_2x_2} Y_t\Big]   - [D_{x_1} Y_t]^+ + 2 [D_{x_2} Y_t]^- =0. 
\eeaa
It is clear that the comparison principle for the viscosity solutions of above PDE holds.
\qed

\begin{rem}
\label{rem-Buckdahn}
{\rm  (i) The above counterexample stays valid when $\a=0$, and thus the game is deterministic. We note that, even in deterministic case, our weak formulation is different from strong formulation. Indeed, the corresponding state process $X^{W, u, v}$ in weak formulation is: 
\beaa
X^{W, 1, u, v}_t = \int_0^t u(s, X^{W,1,u,v}_\cd, X^{W, 2, u, v}_\cd) ds,\q  X^{W, 2, u, v}_t = \int_0^t v(s, X^{W,1,u,v}_\cd, X^{W, 2, u, v}_\cd) ds.
\eeaa
In particular, $X^{W, 2, u, v}$ depends on $u$ as well. Consequently, given $v$, one cannot define $u$ through \reff{Buckdahn2}.

(ii) In this paper the drift coefficient is $b\si$, see \reff{Xweak2}, so the above deterministic example is not covered in our current framework. However, this  assumption is mainly to ensure the wellposedness of the BSDE \reff{BSDEuv} . When $f=0$, one may define the value processes via conditional expectations, instead of $\cY$. Then we may consider $X$ in the form of \reff{Xweak} and all our results, after appropriate modifications,  will still hold true. In particular, the above deterministic game in weak formulation has a value.
\qed}
\end{rem}
 


\end{document}